\documentclass[12pt]{amsart}
\usepackage{amssymb}
\usepackage{amsmath}
\usepackage{mathabx}
\usepackage{amsthm}
\usepackage{amstext}
\usepackage{amsopn}
\usepackage{mathrsfs}
\usepackage{latexsym}
\usepackage{bm}
\DeclareMathAlphabet{\mathpzc}{OT1}{pzc}{m}{it}
\allowdisplaybreaks
\newtheorem{Definition}{Definition}[section]
\newtheorem{Proposition}{Proposition}[section]
\newtheorem{Lemma}{Lemma}[section]
\newtheorem{Theorem}{Theorem}[section]
\newtheorem{Corollary}{Corollary}[section]
\newtheorem{Remark}{Remark}[section]
\newtheorem{Example}{Example}[section]

\begin{document}
%%%%%%%%%%%%%%%%%%%%%%%%%%%%%%%%%%%%%%%%%%%%%%%%%%%%%%%%%%%%
\bibliographystyle{plain}
\footnotetext{
This work is supported by Narodowe Centrum Nauki grant No. 2014/15/B/ST1/00166}
%%%%%%%%%%%%%%%%%%%%%%%%%%%%%%%%%%%%%%%%%%%%%%%%%%%%%%%%%%%%%%%%%
\title[Asymptotic distributions of Wishart type products]
{Asymptotic distributions of Wishart type\\ products of random matrices}
\author[R. Lenczewski, R. Sa{\l}apata]{Romuald Lenczewski and Rafa{\l} Sa{\l}apata}
\address{Romuald Lenczewski \newline
Wydzia\l{} Matematyki, Politechnika Wroc\l{}awska, \newline
Wybrze\.{z}e Wyspia\'{n}skiego 27, 50-370 Wroc{\l}aw, Poland 
\newline
\indent
Rafa{\l} Sa{\l}apata \newline
Wydzia\l{} Matematyki, Politechnika Wroc\l{}awska, \newline
Wybrze\.{z}e Wyspia\'{n}skiego 27, 50-370 Wroc{\l}aw, Poland  \vspace{10pt}}
\email{Romuald.Lenczewski@pwr.edu.pl, Rafal.Salapata@pwr.edu.pl}
%%%%%%%%%%%%%%%%%%%%%%%%%%%%%%%%%%%%%%%%%%%%%%%%%%%%%%%%%%%%%%%%
\begin{abstract}
We study asymptotic distributions of large dimensional random matrices of the form $BB^{*}$, where $B$ 
is a product of $p$ rectangular random matrices, using free probability and combinatorics 
of {\it colored labeled noncrossing partitions}.
These matrices are taken from the set of off-diagonal blocks of the family $\mathpzc{Y}$ 
of independent Hermitian random matrices which are 
asymptotically free, asymptotically free against
the family of deterministic diagonal matrices, and whose norms are uniformly 
bounded almost surely.
This class includes unitarily invariant Hermitian random matrices
with limit distributions given by compactly supported probability 
measures $\nu$ on the real line. 
We express the limit moments 
in terms of colored labeled noncrossing pair partitions, to which we assign 
weights depending on even free cumulants of $\nu$ and on asymptotic dimensions 
of blocks ({\it Gaussianization}). 
For products of $p$ independent blocks, we show that the limit moments 
are linear combinations of a new family of polynomials called {\it generalized multivariate Fuss-Narayana polynomials}.
In turn, the product of two blocks of the same matrix leads to an example with 
rescaled Raney numbers. 
\end{abstract}

\maketitle
\section{Introduction}
The study of products of random matrices has received new attention due to 
their various applications in mathematics, physics and wireless 
communications (some of these applications are listed in [1, 8, 11, 22, 27]).

An important class of products is given by $n\times n$ Wishart matrices $BB^*$ and their product generalizations. Their spectral distributions and asymptotic spectral distributions as $n\rightarrow \infty$ have been studied by many authors 
(see, for instance, [1, 2, 3, 4, 5, 6, 7, 8, 11, 13, 14, 16, 19, 22, 24, 27]). Some of these papers are concerned with products of Gaussian random matrices (GRM, often called the Ginibre Ensemble), as well as products containing not only GRM, but also their inverses or truncated Haar unitary matrices. There is a number of results on random matrices which are more general than GRM, but whose asymptotics is the same as that for GRM (universality).

In this paper, we are going to study asymptotic spectral distributions of 
large dimensional {\it Wishart type products} $BB^*$, where 
$$
B=X_1X_2 \ldots X_p,
$$
and $X_1,X_2, \ldots,  X_p$ are rectangular random matrices of type specified below.
Equivalently, we can speak of the asymptotic distributions of squared singular values of $B$.
All matrices $X_j$ depend on $n\in \mathbb{N}$, which is either equal to 
the first dimension of the first matrix or the dimension 
of some larger square matrix, in which they are embedded as blocks
whose dimensions grow proportionately to $n$. 

Our approach is based on the method of moments and thus we shall investigate the limits
$$
m_{k}:=\lim_{n\rightarrow \infty}\tau((BB^*)^{k})
$$ 
for any $k\in \mathbb{N}$, where 
$\tau$ is the normalized trace composed with classical expectation. 
These moments are interesting objects to study 
from a combinatorial point of view even if the densities of 
the asymptotic distributions are difficult or impossible to compute.
A convenient approach is to embed the matrices in a Hermitian or 
non-Hermitian matrix as blocks lying above the main diagonal 
and realize the expectation of the normalized partial trace (over the subset of 
basis vectors defining the first diagonal block) of the powers of $BB^*$ as
moments of certain operators living in the matricially free Fock space.

In the simplest case, when $B$ is just one $n\times n_1$ matrix with independent 
$N(0,1/n)$ Gaussian entries, the matrix $BB^*$ is the simplest {\it Wishart matrix} and
the limit moments are the well known {\it Narayana polynomials} 
$N_k(t)$ with the associated Marchenko-Pastur distribution of shape $t>0$, 
$$
\rho_{t}={\rm max}\{1-t,0\}\delta_{0}+\frac{\sqrt{(x-a)(b-x)}}{2\pi x}1\!\!1_{[a,b]}(x)dx
$$
where $a=(1-\sqrt{t})^{2}$ and $b=(1+\sqrt{t})^{2}$ and 
$t=\lim_{n\rightarrow\infty}(n_1/n)$, with the understanding that  
$n_1$ depends on $n$. An immediate connection with free probability [28] follows from the 
fact that $\rho_1$ plays the role of the free Poisson law.

This connection becomes much deeper if $B$ is a product of $p$ 
independent $N(0,1/n)$ Gaussian random matrices, $X_1, \ldots , X_p$, 
of sizes $n\times n_1, \ldots , n_{p-1}\times n_p$, respectively, where 
each $n_j$ is proportional to $n$. Then the asymptotic distribution is the
free multiplicative convolution of Marchenko-Pastur laws with various shapes, namely
$$
\rho_{t_1}\boxtimes \rho_{t_2}\boxtimes \ldots \boxtimes \rho_{t_p},
$$
where $t_j=\lim_{n\rightarrow \infty}n_j/n$ for $1\leq j \leq p$.
The moments of this convolution and thus the limit moments of $BB^*$ under the first 
(normalized by $n$) partial trace composed with classical expectation are 
{\it multivariate Fuss-Narayana polynomials}, defined in [16] and studied in the explicit form 
in [19] (see also earlier results in the computer science literature [23]). 

As their name suggests, these polynomials generalize {\it Fuss-Narayana polynomials}, reproduced when all shape parameters are equal, studied earlier in the context of the family of free Bessel laws by Banica, Belinschi, Capitaine and Collins [4]. 
In particular, these authors showed that when all matrices are square, the limit moments are the well known {\it Fuss-Catalan numbers}
$$
F_{k}(p) = 
\frac{1}{k}{(p+1)k \choose pk+1},
$$
obtained also in the study of powers and products of 
a larger than Gaussian class of independent random matrices by Alexeev, G\"{o}tze and Tikhomirov [2, 3].
These numbers define probability distributions whose 
densities were computed in the explicit form by Penson and $\dot{\rm Z}$yczkowski [27]. 

In the case of independent GRM, we can take matrices $X_j$ to be
certain off-diagonal blocks of one large hermitian GRM $Y$ (with i.i.d. or i.b.i.d. entries, where i.b.i.d. stands for independent, block identically distributed). Observe that the same result is obtained when each $X_j$ is a block 
of a different matrix taken from a family of independent GRM.
In the more general case, when the blocks of each $Y$ are not necessarily independent, 
this is no longer the case. In  particular, it is natural to take $X_j$ to be the blocks lying above the main diagonal of 
a quite general family of independent random matrices, with squared singular value distributions much more general than Marchenko-Pastur
laws. In this context, the blocks are independent if they are taken from different matrices, 
but may be dependent when they are taken from the same matrix.

Therefore, in this paper, we replace the family of independent blocks of GRM by 
off-diagonal blocks of a quite general family of independent random matrices,
$$
\mathpzc{Y}=\{Y(u): u\in \mathpzc{U}\},
$$ 
which is asymptotically free, asymptotically free against the family of constant 
diagonal matrices under $\tau$ and whose norms are uniformly bounded almost surely. 
This class includes unitarily invariant random matrices whose 
asymptotic spectral distributions are given by compactly supported probability measures $\nu(u)$ on the real line
(in particular, the Gaussian Unitary Ensemble and Wigner matrices). In general, the blocks of one such matrix 
may be dependent.

We give a unified operatorial and combinatorial 
description of the limit moments of the class of Wishart type products constructed
from blocks of $\mathpzc{Y}$.
The operatorial description consists in expressing asymptotic distributions 
in terms of mixed moments of variables that can be viewed as
matricial counterparts of canonical noncommutative random variables. 
The combinatorial one consists in expressing 
the limit moments in terms of {\it colored labeled noncrossing pair partitions}, to which
we assign weights depending on even free cumulants of 
$\nu(u)$, where $u$ is a {\it label}, and on asymptotic dimensions $d_j$ 
of blocks, where $j$ is a {\it color}.
This result is based on a bijection between the set of {\it all} noncrossing 
partitions {\it adapted} to words $W^k$ that encode information about 
$m_k$ and noncrossing {\it pair} partitions adapted to
twice longer words $\widetilde{W}^k$ and plays the role of 
a {\it Gaussianization} of the asymptotic distributions. 
In our notation,
$$
\mathcal{NC}(W^{k})\cong \mathcal{NC}^{2}(\widetilde{W}^{k}),
$$
which is natural since we have the well-known bijection $\mathcal{NC}(n)\cong \mathcal{NC}^{2}(2n)$, 
but it is nontrivial due to the additional content related to adaptedeness. 
In this scheme, we assign free cumulants to 
blocks of odd depths and suitable asymptotic dimensions to blocks of even depths.

These combinatorial results hold in the general case, when blocks 
are taken from the same matrix or from different matrices.
Then we study two situations in more detail.
First, we assume that the product matrix $B$ is built from off-diagonal
blocks of different independent matrices that belong to $\mathpzc{Y}$ (for any $p$), generalizing 
the case of blocks of (one or many, in the Gaussian case this is not relevant) 
independent GRM. 
In this context, we find a formula for the moment generating function 
for the limit moments and we show that they are linear combinations of 
multivariate polynomials of the form 
$$
P_{m,r}(d_1,d_2,\ldots,d_{p+1}) = \sum_{j_1+\ldots+j_{p+1}=mp+r}
 \frac1n 
{m \choose j_1}\ldots {m\choose j_{p+1}}
d_1^{j_1} \ldots d_{p+1}^{j_{p+1}}.
$$
where $m,r\in {\mathbb N}$ and summation runs over nonnegative integers, called 
{\it generalized multivariate Fuss-Narayana polynomials} (for $r=1$ we obtain 
multivariate Fuss-Narayana polynomials). The limit moments are shown to be linear combinations 
of these polynomials. One can say that the presence of polynomials with $r\neq 1$ means that 
the limit moments differ from those for GRM.

Second, we assume that $B$ is a product of blocks 
of the same matrix $Y$ and thus may be highly dependent. This case is studied for $p=2$ 
(analytic and numerical results for larger $p$ are technically
more involved). A formula for the moment generating function is derived.  
We also show that if the asymptotic distribution of $Y$ is the Marchenko-Pastur distribution
and all asymptotic dimensions are equal to one, then 
the limit moments $m_k$ are certain rescaled Raney numbers. By construction, they are moments of a probability measure.
The general conditions under which given Raney numbers define a probability measure were given by M{\l}otkowski 
[21]. The corresponding densities have also been computed, with the most general result given by Forrester and Liu [12].

We use operatorial realizations of the limit moments in terms of
matricial counterparts of canonical noncommutative random variables of Voiculescu.
The formalism of matricial freeness, used in this paper, can be viewed as freeness with respect to 
a family of scalar-valued states and was developed by one of the authors 
and applied to asymptotic distributions of blocks of random matrices in [16, 17, 18].
This gave an alternative approach to that of Benaych-Georges in [6,7], who earlier described 
the collective behavior of independent rectangular random matrices whose sizes grow to infinity at different
rates. His approach was based on the formalism of freeness with amalgamation over the algebra $A$ of diagonal matrices
and the associated conditional expectation and developed some analytic tools (like the rectangular $R$-transform) to study 
the asymptotic distributions of their sums. As shown in [18], the two formalisms are related 
to each other. In particular, canonical operators on the matricially free Fock space 
can be expressed in terms of the canonical operators 
on the full Fock space over a suitably defined Hilbert $A$-bimodule. 

Besides this introduction, this paper contains 6 more sections.
Section 2 contains information about operatorial tools. 
Section 3 describes limit moments of matrices of Wishart type,
using the combinatorics of colored labeled noncrossing partitions adapted to 
words $W$.
Section 4 discusses the Gaussianization of asymptotic distributions, by which 
the limit moments are described in terms of colored labeled noncrossing pair partitions 
adapted to twice longer words $\widetilde{W}$. Section 5 treats the limit moments of products of 
independent blocks and expresses them in terms of generalized multivariate Fuss-Narayana polynomials.
Section 6 contains two examples. Section 7 treats the limit moments of products of 
two (in general, dependent) blocks of the same matrix.

For any natural $n$, we use the notation $[n]=\{1,2, \ldots , n\}$. 
We also assume that $p\in {\mathbb N}$, often suppressed in our notations.

\section{Operators}
Let us study the case when $B$ is the product of
rectangular random matrices taken to be the off-diagonal 
blocks taken from the family of independent random matrices, about which 
we make some assumptions. 
We would like to study blocks of the family of matrices, to which we can apply our previous results on the asymptotics of blocks.  

In order to apply the framework of free probability to asymptotic 
random matrix theory, we usually consider the normalized trace 
composed with classical expectation, namely 
$$
\tau(A)=\frac{1}{n}{\mathbb E}({\rm Tr}(A)),
$$
for any $A\in M_{n}({\mathbb C})$. 

We will consider the family $\mathpzc{Y}=\{Y(u):u\in \mathpzc{U}\}$ 
of independent Hermitian random matrices and we will assume that 
\begin{enumerate}
\item[(A1)]
$\mathpzc{Y}$ is asymptotically free under $\tau$,
\item[(A2)]
$\mathpzc{Y}$ is asymptotically free under $\tau$ 
against the family of deterministic diagonal matrices,
\item[(A3)]
the norms of $Y(u)\in \mathpzc{Y}$ are uniformly bounded almost surely.
\end{enumerate}
These three assumptions will be made in the sequel and will not always be explicitly stated. 
The class of matrices under consideration will include, in particular, 
unitarily invariant Hermitian matrices with asymptotic distributions given by 
compactly supported probability measures $\nu(u)$ on the real line, as well as some other classical ensembles 
like Wigner matrices. Such objects as 
matrices, their blocks and traces will always depend on $n$, but this 
will be suppressed in our notations for simplicity. 

In the framework of free probability, if we have such a family 
$\{Y(u):u\in \mathpzc{U}\}$, we will symbolically write  
$$
\lim_{n\rightarrow \infty}Y(u)=\gamma(u),
$$
where convergence is in the sense of mixed moments under $\tau$ and each 
$\gamma(u)$ is a {\it canonical noncommutative random variable} 
for any $u\in \mathpzc{U}$.
This variable can be written in the form
$$
\gamma(u)=\wp(u)^*+\sum_{k=0}^{\infty}r_{k+1}(u)\wp(u)^{k},
$$
where $\{\wp(u):u\in \mathpzc{U}\}$ is a family of free creation operators and 
$(r_{k}(u))$ is the sequence of free cumulants of the limit measure of $Y(u)$
for any $u$. In general, this series can be treated as a formal series, living in a 
suitably defined noncommutative probability space [28], but it is a 
well-defined bounded operator on the free Fock space
if the series of free cumulants is absolutely convergent.

We will consider the matrices of Wishart type $B$ built from
off-diagonal blocks of matrices $Y(u)$. Therefore, we denote 
$$
S_{i,j}(u)=D_{i}Y(u)D_{j}
$$
for any $i,j\in [r]$, where $I_n=D_{1}+D_2+\ldots +D_{r}$ 
is the natural partition of the $n\times n$ unit matrix, namely the only nonvanishing 
entries are $(D_j)_{k,k}=1$ for $n_{1}+\ldots +n_{j-1}< k\leq n_1+\ldots + n_j$, where 
$n_1+\ldots +n_{r}=n$. Again, as in most other cases, the dependence on $n$ of matrices
$D_j$ and numbers $n_j$ is suppressed. 

In the study of products of random matrices it is convenient to use 
blocks lying above the main diagonal and assume that $r=p+1$ for natural $p$.
In that setting, by {\it asymptotic dimensions} we shall understand 
numbers
$$
d_{j}=\lim_{n\rightarrow \infty}\frac{n_{j}}{n}
$$
for $1\leq j \leq p+1$. Thus, we get the matrix $D={\rm diag}(d_1, \ldots , d_{p+1})$ of trace one.
In the context of multivariate Fuss-Narayana polynomials [19] it was slightly more convenient to use $d_0, d_1, \ldots , d_p$
and renormalize all moments by setting $n=n_{0}$, which gives $d_0=1$, the remaining asymptotic dimensions 
being arbitrary positive numbers.

We described in [16] the limit joint distributions of such blocks
in terms of mixed moments of a family
of bounded operators on the {\it matricially free Fock space of tracial type}
${\mathcal M}$. 
This Fock space is a direct sum Hilbert space built from arrays of Hilbert spaces
$\{{\mathcal H}_{i,j}(u):(i,j)\in \mathpzc{J}, u\in\mathpzc{U}\}$, where $\mathpzc{J}\subset [r]\times [r]$, of the form
\begin{equation*}
{\mathcal M}= \bigoplus_{j=1}^{r} {\mathcal M}_{j},
\end{equation*}
where 
\begin{equation*}
{\mathcal M}_{j}={\mathbb C}\Omega_{j}\oplus \bigoplus_{m=1}^{\infty}
\bigoplus_{\stackrel{i_1,\ldots , i_m}
{\scriptscriptstyle u_1, \ldots , u_n}
}
{\mathcal H}_{i_1,i_2}(u_{1})\otimes {\mathcal H}_{i_2,i_3}(u_2)
\otimes \ldots \otimes 
{\mathcal H}_{i_m,j}(u_{m}),
\end{equation*}
for any $j$, with $\Omega_j$ being a unit vector, endowed with the canonical inner product.
Denote by $P_j$ the canonical orthogonal projection onto ${\mathbb C}\Omega_j$.
In this work, it suffices to assume that ${\mathcal H}_{i,j}(u)={\mathbb C}e_{i,j}(u)$, where 
$\{e_{i,j}(u):i,j \in [r], u\in \mathpzc{U}\}$ is an orthonormal set.
Note that the space $\mathcal{M}$ is a simple version of the `free Fock space of matricial 
type' needed in the random matrix context.

We found a realization of $\gamma(u)$ as an operator on 
${\mathcal M}$ and showed in [16, 18] that $D_j$ tends to $P_j$ for any $j$.
Thus, we can write
$$
\lim_{n\rightarrow \infty}S_{i,j}(u)=P_i\gamma(u)P_j
$$
for any $i,j,u$, where convergence is in mixed moments under $\tau$. 
In fact, when speaking of mixed moments of blocks, it is natural to 
use partial traces $\tau_j$ over the sets of $n_j$ basis vectors. Then
we have to consider the family of scalar-valued states $\{\Psi_j: 1\leq j \leq r\}$ on 
$\mathcal{M}$, where 
$$
\Psi_j(.)=\langle .\Omega_j, \Omega_j\rangle
$$ 
for any $j$ and the asymptotic distributions of blocks under $\tau_j$ are described 
by the mixed moments of operators on ${\mathcal M}$ under $\Psi_j$. 

In this approach, the basic operators are $\wp_{i,j}(u)=P_i\wp(u)P_j$, called 
the {\it matricially free creation operators}. Their non-trivial action onto the basis vectors of $\mathcal{M}$
is given by
\begin{eqnarray*}
\wp_{i,j}(u)\Omega_j&=&\sqrt{d_{i}}e_{i,j}(u)\\
\wp_{i,j}(u)(e_{j,k}(s))&=&\sqrt{d_{i}}(e_{i,j}(u)\otimes e_{j,k}(s))\\
\wp_{i,j}(u)(e_{j,k}(s)\otimes w)&=&\sqrt{d_{i}}(e_{i,j}(u)\otimes e_{j,k}(s)\otimes w)
\end{eqnarray*}
for any $i,j,k\in [r]$ and $u,s\in \mathpzc{U}$, where $e_{j,k}(s) \otimes w$ is assumed to be 
a basis vector, with the understanding that their actions onto the remaining basis vectors give zero. 
A similar definition, with a more general scaling factor, was used in [16, 18].
The corresponding {\it matricially free annihilation operators} are their adjoints,
denoted $\wp_{p,q}(u)^{*}$. 

Let us remark that there is an alternative operator-valued
approach in terms of conditional expectation and Hilbert bimodules. 
In the operator-valued free probability, the limit mixed moments of blocks under $\tau$
are expressed in terms of mixed moments under conditional expectation of 
the canonical isometries $\ell(u)$ on the free Fock space ${\mathcal F}({\mathcal H})$ over a 
Hilbert $A$-bimodule ${\mathcal H}$, where $A$ is the $C^*$-algebra of diagonal $r\times r$ matrices, and of the 
canonical generators $F_1, \ldots, F_r$ of $A$.
Roughly speaking, our operators $\wp_{i,j}(u)$ correspond then 
to operators of the form $F_i\ell(u)F_j$, where $i,j\in [r]$ and $u\in \mathpzc{U}$. 
More details on the correspondence between matricial freeness and 
freeness with amalgamation over $A$ can be found in [18]. 

\begin{Definition}
{\rm Let $(r_{k}(u))$ be the sequence of free cumulants of the distribution of $\gamma(u)$.
Let us introduce operators $\Gamma_{i,j}(u)=P_i\gamma(u)P_j$, which take 
the form
\begin{eqnarray*}
\Gamma_{i,j}(u)&=&\wp_{j,i}(u)^{*}+\delta_{i,j}r_{1}(u)P_{j}\\
&&+\sum_{k=1}^{\infty}r_{k+1}(u)
\sum_{i_1, \ldots , i_{k-1}}\wp_{i,i_1}(u)\wp_{i_1,i_2}(u)\ldots \wp_{i_k,j}(u)
\end{eqnarray*}
where $i,j\in [r]$ and $u\in \mathpzc{U}$. Moreover, we set 
$\Gamma_{i,j}^*(u)=\Gamma_{j,i}(u)$ for any $i,j,u$ (this is just a convenient notation). }
\end{Definition}

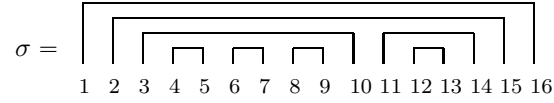
\begin{figure}
\unitlength=1mm
\special{em:linewidth 1pt}
\linethickness{0.4pt}
\begin{picture}(70.00,20.00)(-8.00,0.00)
%%%%%%%%%%%%%%%%%%%%%%%%%%%%%%%%%%%%%%%%%%%%%%%%%%%%%%%%555
\put(-4.00,12.00){\line(1,0){60.00}}
\put(0.00,10.00){\line(1,0){52.00}}
\put(4.00,8.00){\line(1,0){28.00}}
\put(8.00,6.00){\line(1,0){4.00}}
\put(16.00,6.00){\line(1,0){4.00}}
\put(24.00,6.00){\line(1,0){4.00}}
\put(40.00,6.00){\line(1,0){4.00}}
\put(36.00,8.00){\line(1,0){12.00}}
%%%%%%%%%%%%%%%%%%%%%%%%%%%%%%%%%%%%%%%%%%%%%%%%%%%%%%%%%%%
\put(-4.00,4.00){\line(0,1){8.00}}
\put(0.00,4.00){\line(0,1){6.00}}
\put(4.00,4.00){\line(0,1){4.00}}
\put(8.00,4.00){\line(0,1){2.00}}
\put(12.00,4.00){\line(0,1){2.00}}
\put(16.00,4.00){\line(0,1){2.00}}
\put(20.00,4.00){\line(0,1){2.00}}
\put(24.00,4.00){\line(0,1){2.00}}
\put(28.00,4.00){\line(0,1){2.00}}
\put(32.00,4.00){\line(0,1){4.00}}
\put(36.00,4.00){\line(0,1){4.00}}
\put(40.00,4.00){\line(0,1){2.00}}
\put(44.00,4.00){\line(0,1){2.00}}
\put(48.00,4.00){\line(0,1){4.00}}
\put(52.00,4.00){\line(0,1){6.00}}
\put(56.00,4.00){\line(0,1){8.00}}

\put(-13.00,5.00){\footnotesize $\sigma=$}
\put(-4.50,0.00){$\scriptstyle{1}$}
\put(-0.50,0.00){$\scriptstyle{2}$}
\put(3.50,0.00){$\scriptstyle{3}$}
\put(7.50,0.00){$\scriptstyle{4}$}
\put(11.50,0.00){$\scriptstyle{5}$}
\put(15.50,0.00){$\scriptstyle{6}$}
\put(19.50,0.00){$\scriptstyle{7}$}
\put(23.50,0.00){$\scriptstyle{8}$}
\put(27.50,0.00){$\scriptstyle{9}$}
\put(31.50,0.00){$\scriptstyle{10}$}
\put(35.50,0.00){$\scriptstyle{11}$}
\put(39.50,0.00){$\scriptstyle{12}$}
\put(43.50,0.00){$\scriptstyle{13}$}
\put(47.50,0.00){$\scriptstyle{14}$}
\put(51.50,0.00){$\scriptstyle{15}$}
\put(55.50,0.00){$\scriptstyle{16}$}
\end{picture}
\caption{A noncrossing pair partition of Example 3.1. }
\end{figure}

\begin{Remark}
{\rm The operators $\Gamma_{i,j}(u)$ play the role of matricial counterparts of 
canonical noncommutative random variables. They are similar to slightly different operators 
\begin{eqnarray*}
\gamma_{i,j}(u)&=&\wp_{i,j}(u)^{*}+\delta_{i,j}r_{1}(u)P_{j}\\
&&+\sum_{k=1}^{\infty}r_{k+1}(u)
\sum_{i_1, \ldots , i_{k-1}}
\wp_{i,i_1}(u)\wp_{i_1,i_2}(u)\ldots \wp_{i_k,j}(u),
\end{eqnarray*}
introduced in [16], with which they coincide for $i=j$.}
\end{Remark}

\section {Combinatorics of limit moments}

Let us recall certain basic definitions. If we are given a noncrossing partition $\pi$ and $V,V'$ 
are blocks of $\pi$, then $V'$ is an {\it inner} block of $V$ if $V=A\cup B$ for some disjoint nonempty subsets $A$ and $B$,
and it holds that $i<j<k$ for any $i\in A, k\in B$ and $j\in V'$ (then $V$ is an {\it outer} block of $V'$).
Moreover, $V'$ is the {\it nearest inner block} of $V$ if it is an inner block and there does not exist
a block $V''$ such that $V'$ is an inner block of $V''$ and $V''$ is an inner block of $V$.  
We denote by ${\mathpzc i}(V)$ the number of nearest inner blocks of $V$. In a similar way we can define 
the nearest outer block of a given block. We will say that 
$V$ has depth $\mathpzc{d}$ if its nearest outer block has depth $\mathpzc{d}-1$, assuming 
that a block without outer blocks has depth $1$.
The {\it depth} of $V$ will be denoted $\mathpzc{d}(V)$.

\begin{Example}
{\rm The noncrossing pair partition $\sigma$ given in Fig. 1 consists of 8 blocks, among which there are 3 blocks of odd depths and 
5 blocks of even depths. For instance, $\mathpzc{d}(\{1,16\})=1$, $\mathpzc{d}(\{2,15\})=2$, $\mathpzc{d}(\{3,10\})=3$
and $\mathpzc{d}(\{4,5\})=4$. Blocks $\{4,5\}, \{6,7\},\{8,9\}$ are the nearest inner blocks of $\{3,10\}$ 
and $\{12,13\}$ is the nearest inner block of $\{11,14\}$. Since numbers $\mathpzc{i}(V)$ will be important 
for blocks of odd depths, we compute 
$\mathpzc{i}(\{1,16\})=\mathpzc{i}(\{11,14\})=1$ and 
$\mathpzc{i}(\{3,10\})=3$. Let us add that the diagram in Fig. 1 and the right diagram in Fig. 5 represent
the same partition, except that in the circular version in Fig. 5 we speak of vertices and edges instead of
legs and segments, respectively. We adopt the convention that the bottom vertex of label $1$ and color $1$, 
namely $(16,1)$, corresponds to leg $1$ in Fig. 1 and the vertices corresponding to the 
remaining legs are placed counterclockwise.}
\end{Example}

The combinatorics of limit moments of blocks is based on {\it colored labeled noncrossing partitions}, by which we understand noncrossing partitions whose legs are {\it colored} (or, *-colored) by a certain *-alphabet $A_p$ 
and equipped with {\it labels} from the set $\mathpzc{U}$ carrying information about 
independence of the corresponding matrices [16]. Let us adapt these results to the case of matrices of Wishart type.
\begin{enumerate}
\item[(a)]
For fixed $p$, the $k$th moment of $BB^*$ is a mixed moment of $2pk$ matrices. Therefore, it is 
natural to expect that the limit moment will be expressed in terms of noncrossing partitions 
$\pi$ of the set of $2pk$ numbers. We will use the alphabet 
$$
A_{p}=\{(j,u): j\in \{1, \ldots , p,1,\ldots , p^*\}, \;u\in \mathpzc{U}\},
$$
where $\mathpzc{U}$ is a finite set, to encode some additional information about these partitions. The letters of these alphabets will be assigned to the legs of $\pi$ and we will call $j$ its {\it color}, 
whereas $u$ will be its {\it label}. Different labels correspond to different independent matrices, but we assume that 
the same label $u$ is associated with $j$ and $j^*$. The same labels indicate that matrices are blocks of the same matrix.
\item[(b)]
For simplicity, we will usually write $j$ instead of $(j,u)$. Thus, let us introduce the basic word
$$
W=12\ldots pp^*\ldots 2^*1^*,
$$
then we will say that $\pi\in \mathcal{NC}(2pk)$ is {\it label adapted} to the word $W^k$ if and only if 
the tuple assigned to each block $V$ is of the form $((j_1,u), \ldots ,(j_m,u))$, where $j_1, \ldots , j_m\in A_{p}$
and $u\in \mathpzc{U}$. In this case we will write $u=u_{V}$.
\begin{figure}
\unitlength=1mm
\special{em:linewidth 1pt}
\linethickness{0.4pt}
\begin{picture}(55.00,24.00)(-05.00,00.00)
%%%%%%%%%%%%%%%%%%%%%%%%%%%%%%%%%%%%%%%%%%%%%%%%%%%%%%%%555
\put(0.00,10.00){\line(1,0){50.00}}

\put(0.00,5.00){\line(0,1){5.00}}
\put(10.00,5.00){\line(0,1){5.00}}
\put(20.00,5.00){\line(0,1){5.00}}
\put(40.00,5.00){\line(0,1){5.00}}
\put(50.00,5.00){\line(0,1){5.00}}

\put(-10.00,5.00){\footnotesize $V=$}
\put(-01.00,1.00){\footnotesize $j_1$}
\put(09.00,1.00){\footnotesize $j_2$}
\put(19.00,1.00){\footnotesize $j_3$}
\put(37.00,1.00){\footnotesize $j_{m-1}$}
\put(49.00,1.00){\footnotesize $j_m$}
\put(-10.00,5.00){\footnotesize $V=$}
\put(04.00,12.00){\footnotesize $i_1$}
\put(14.00,12.00){\footnotesize $i_2$}
\put(29.00,5.00){\footnotesize $\ldots$}
\put(24.00,12.00){\footnotesize $i_3\;\;\ldots$}
\put(42.00,12.00){\footnotesize $i_{m-1}$}
\end{picture}
\caption{Block with colored legs and segments}
\end{figure}
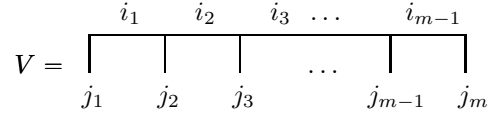

\item[(c)]
As in [16], we not only assign colors to the legs of $V$, but we also color 
its subblocks consisting of two consecutive numbers. These will be represented by horizontal 
segments connecting neighboring legs, as shown in Fig. 2. In this fashion we obtain a multicolored block 
$V$ with coloring of segments given by the formula
$$
i_m=\left\{
\begin{array}{ll}
j+1 & {\rm if}\; j_m=j\\
j   & {\rm if}\; j_m=j^*
\end{array}
\right.
$$
which allows us to assign to $V$ the {\it coloring} in the form of the word
$$
c(V)=i_1i_2\ldots i_{m-1}.
$$
By the above definition, the set of colors used for segments is slightly different than that for legs,
where we also use colors, namely $C=\{1,\ldots , p+1\}$ for given $p$.
\item[(d)]
Coloring segments is slightly more convenient than coloring legs of $V$ if we want to
encode information about the moments of the corresponding operators. This is because, as we showed in [16],
the product of operators
$$
\wp_{i_1,j}^*\wp_{i_2,i_1}^*\ldots \wp_{i_{m-1},i_{m-2}}^{*}\wp_{i_{m-1},i_{m-2}}\ldots \wp_{i_2,i_1}\wp_{i_{1},j}
$$
is naturally associated with an $m$-block $V$ with coloring 
$i_1i_2\ldots i_{m-1}$. The first $m-1$ legs correspond to the annihilation operators, 
and the rightmost leg corresponds to the product of creation operators, as shown in Fig. 3.
The above product contributes the product of asymptotic dimensions
$$
d(V)=d_{i_1}d_{i_2}\ldots d_{i_{m-1}},
$$
called the {\it dimension weight}, defined by the coloring $c(V)$.
\item[(e)]
An important property of the above product of operators is that in the Hilbert space 
framework it acts non-trivially onto vectors from ${\mathcal M}$ with first index 
equal to $j$, namely $e_{j,l}(u)$ or $e_{j,l}(u)\otimes w$ for some $w$. 
In segment coloring, $j$ is then assigned to the segment lying above block $V$ and 
belonging to the nearest outer block of $V$. 
\item[(f)]
Recall that any noncrossing partition of $[n]$ can be defined by bracketing the sequence 
of numbers $12\ldots n$, identified with ${\bf 1}_{n}$, the coarsest partition of $[n]$
consisting of one block. This bracketing, called noncrossing, must be such that the 
pairs of brackets (left and right) form a noncrossing pairing. 
We will apply the same kind of bracketing with some additional conditions 
in order to define a suitable class of noncrossing partitions of the set $[2pk]$ that 
is adapted to the colorings of their blocks.
\item[(g)]
Consider the coarsest partition ${\bf 1}_{2pk}$ of $[2pk]$ consisting of one block, 
for any given $m$ and $p$. Its coloring is very simple, namely 
$$
c({\bf 1}_{2pk})=\underbrace{c_p1c_p1\ldots c_p1c_p}_{c_p \;\;appears\;k \;\;times}
$$
where $c_p=2\ldots p(p+1)p\ldots 2$. By adding $1$ both at the beginning and the end of the coloring, we 
obtain the sequence 
$$
\widetilde{c}({\bf 1}_{2pk})=\underbrace{1c_p1c_p1\ldots c_p1c_p1}_{c_p \;\;appears\;k \;\;times}
$$
that can be identified with the unique periodic Catalan path of length $2pk$ and period $2p$, given by 
the real-valued function $f:[0,2pk]\rightarrow [0,p]$ with the usual north-east and south-east steps, 
for which  $f(j)=i_{j-1}-1$. 
\item[(h)]
A noncrossing bracketing of the sequence $\widetilde{c}({\bf 1}_{2pk})$, defined analogously to 
a noncrossing bracketing of $[n]$, will be called {\it color adapted}
if and only if each pair of brackets (in fact, we use parentheses rather than brackets) 
has identical neighbors, by which we mean that locally it has the form 
$$
i(i_1\ldots i_{m-1})i
$$
for some $i$, where under $...$ we may also have brackets with some sequence of numbers. 
In particular, $m$ must be even since $|i_s-i_{s+1}|=1$ for $1\leq s\leq m-2$ 
and $|i_1-i_{m-1}|\in \{0,2\}$. This implies that by applying a sequence of operations 
to the sequence with a noncrossing bracketing, consisting in removing the expressions in brackets which do not have any brackets 
inside and identifying the neighboring colors $i$, we can reduce 
the coloring $\widetilde{c}({\bf 1}_{2pk})$ to the coloring $1$ (see Example 3.3).
Note that the first step of this reduction corresponds to 
removing the product of operators considered in (d), keeping in mind that 
this product contributes the numerical factor $d(V)$.
\begin{figure}
\unitlength=1mm
\special{em:linewidth 1pt}
\linethickness{0.4pt}
\begin{picture}(55.00,20.00)(0.00,-4.00)
%%%%%%%%%%%%%%%%%%%%%%%%%%%%%%%%%%%%%%%%%%%%%%%%%%%%%%%%555
\put(0.00,10.00){\line(1,0){60.00}}

\put(0.00,5.00){\line(0,1){5.00}}
\put(10.00,5.00){\line(0,1){5.00}}
\put(20.00,5.00){\line(0,1){5.00}}
\put(40.00,5.00){\line(0,1){5.00}}
\put(60.00,5.00){\line(0,1){5.00}}

\put(-10.00,5.00){\footnotesize $V=$}
\put(-03.00,1.00){\footnotesize $\wp_{i_1,j}^*$}
\put(07.00,1.00){\footnotesize $\wp_{i_2,i_1}^*$}
\put(17.00,1.00){\footnotesize $\wp_{i_3,i_2}^{*}$}
\put(28.00,5.00){\footnotesize $\ldots$}
\put(38.00,1.00){\footnotesize $\wp_{i_{m-1},i_{m-2}}^{*}$}
\put(56.00,1.00){\footnotesize $\wp_{i_{m-1},i_{m-2}}^{}\ldots \wp_{i_{1},j}^{}$}

\put(04.00,12.00){\footnotesize $i_1$}
\put(14.00,12.00){\footnotesize $i_2$}
\put(24.00,12.00){\footnotesize $i_3\;\;\ldots$}
\put(45.00,12.00){\footnotesize $i_{m-1}$}

\end{picture}
\caption{A product of creation operators is assigned to the rightmost leg and 
annihilation operators are assigned to the remaining legs. }
\end{figure}

\item[(i)]
A partition $\pi\in \mathcal{NC}(2pk)$ will be called {\it color adapted} to $W^k$
if the associated noncrossing bracketing of $\widetilde{c}({\bf 1}_{2pk})$ is color adapted.
It will be called {\it adapted} to the word $W^k$ if and only it is color adapted  
and label adapted to $W^{k}$. This family of partitions 
will be denoted by $\mathcal{NC}(W^k)$. Note that all blocks of such partitions are even.
\item[(j)]
Since each partition of $[n]$ is a collection of blocks and each block, say $V=\{s_1<s_2<\ldots <s_m\}$,
can be represented as the cycle $(s_1,s_2,\ldots ,s_m)$, each partition $\pi\in \mathcal{NC}(n)$ 
can be represented as a product of cycles of this form. In that case, we can treat $\pi$ as a permutation of 
$[n]$, for which $\pi(s_1)=s_2,\ldots , \pi(s_m)=s_1$ for each block $V$ of the above form. 
It can be seen that $\pi\in \mathcal{NC}(2pk)$ is color adapted to $W^k$
if and only if each of its cycles is even and either it is a {\it shift} modulo $2p$, or a {\it reflection} modulo $2p$, 
namely
$$
\pi(s_{j})=s_{j}+1\;\;{\rm mod}(2p),
$$
or
$$
\pi(s_{j})=2p-s_j+1\;\;{\rm mod}(2p),
$$
respectively, for any $j$. In certain cases, both a shift and a reflection 
can be used to produce $\pi(s)$ from $s$. 
\item[(k)]
Let us also observe that the reduction of colorings 
described in (h) can also be carried out on the level of Catalan paths. One starts with the 
periodic Catalan path of length $2pk$ and cuts out the part corresponding to
the bracketed subsequence (in our inductive procedure, that one corresponds to blocks which do not have inner blocks). 
In this fashion one obtains a sequence of Catalan paths which reduces to a point after a finite number of steps.
\end{enumerate}

\begin{figure}
\unitlength=1mm
\special{em.linewidth 1pt}
\linethickness{0.5pt}
\begin{picture}(50.00,55.00)(-85.00,-25.00)
%%%%%%%%%%%%%%%%%%%%%%%%%%%%%%%%%%%%%%%%%%%%
\put(-56.00,-12.00){\scriptsize $d_2$}
\put(-73.00,-06.00){\scriptsize $d_{2}$}
\put(-66.00,-12.00){\scriptsize $r_{8}$}
\put(-50.00,-06.00){\scriptsize $d_3$}
\put(-56.00,11.00){\scriptsize $d_{1}$}
\put(-49.00,04.00){\scriptsize $d_{2}$}
\put(-65.00,11.00){\scriptsize $d_2$}
\put(-73.00,04.00){\scriptsize $d_{3}$}

%%%%%% quadrant 1  circle 1 %%%%%%%%%%%%%%%%%%%%%%%%%%%%%
\put(-60.00,20.00){\circle*{1.00}}
\put(-60.50,21.50){$\scriptstyle{1}$}
\put(-40.00,00.00){\circle*{1.00}}
\put(-38.50,-00.50){$\scriptstyle{2^*}$}
\put(-45.80,14.20){\circle*{1.00}}
\put(-44.80,15.20){$\scriptstyle{1^*}$}

%%%%%%%%%%%%%%%%% quadrant 4 circle 1 %%%%%%%%%%%%
\put(-45.80,-14.20){\circle*{1.00}}
\put(-44.80,-16.00){$\scriptstyle{2}$}
\put(-60.00,-20.00){\circle*{1.00}}
\put(-60.50,-23.00){$\scriptstyle{1}$}

%%%%%%%%%%%%%%%%% quadrant 3 circle 1 %%%%%%%%%%%%
\put(-80.00,00.00){\circle*{1.00}}
\put(-84.00,-00.50){$\scriptstyle{2^*}$}
\put(-74.14,-14.20){\circle*{1.00}}
\put(-77.50,-16.60){$\scriptstyle{1^*}$}

%%%%%%%%%%%%%%%%% quandrant 2 circle 1 %%%%%%%%%%%
\put(-74.10,14.10){\circle*{1.00}}
\put(-76.00,15.50){$\scriptstyle{2}$}

%%%%%%%%%%%%%%% edges circle 1 %%%%%%%%%%%%%%%%%%%%%%%%%%%%
%%%%%%%%%%%%% 1-2 %%%%%%%%%%%%%%%%%%%%%%%%%%%
\qbezier(-60.00,20.00)(-66.00,10.00)(-74.10,14.10)
%%%%%%%%%%%%% 2-2* %%%%%%%%%%%%%%%%%%%%%%%%%%%
\qbezier(-80.00,0.00)(-70.00,07.00)(-74.10,14.10)
%%%%%%%%%%% 2*-1* %%%%%%%%%%%%%%%%%%%%%%%%%%%%%%
\qbezier(-40.00,0.00)(-50.00,07.00)(-45.80,14.20)
%%%%%%%%%%% 1*-1 %%%%%%%%%%%%%%%%%%%%%%%%%%%%%%%
\qbezier(-45.80,14.20)(-55.00,10.00)(-60.00,20.00)
%%%%%%%%%%%% 2*-1* %%%%%%%%%%%%%%%%%%%%%%%%%%%%%%%%%%%%%%%
\qbezier(-80.00,0.00)(-70.00,-7.00)(-74.14,-14.20)
%%%%%%%%%% 1-2 %%%%%%%%%%%%%%%%%%%%%%%%%%%%%%%%%
\qbezier(-60.00,-20.00)(-55.00,-10.00)(-45.80,-14.20)
%%%%%%%%%%% 2-2* %%%%%%%%%%%%%%%%%%%%%%%%%%%%%%%%%%%%%%%
\qbezier(-40.00,0.00)(-50.00,-5.00)(-45.80,-14.20)
%%%%%%%%%%%%% 1*-1 %%%%%%%%%%%%%%%%%%%%%%%%%%%
\qbezier(-60.00,-20.00)(-66.00,-10.00)(-74.14,-14.20)

%%%%%%%%%%%%%%%%%%% circle 1 from bezier curves %%%%%%%%%%%%%%%%%%%%%%%
\qbezier(-60,20)(-52,20)(-46,14)
\qbezier(-46,14)(-40,08)(-40,00)
\qbezier(-40,00)(-40,-08)(-46,-14)
\qbezier(-46,-14)(-52,-20)(-60,-20)
\qbezier(-60,-20)(-68,-20)(-74,-14)
\qbezier(-74,-14)(-80,-08)(-80,00)
\qbezier(-80,00)(-80,08)(-74,14)
\qbezier(-74,14)(-68,20)(-60,20)
\end{picture}
\caption{Cycle corresponding to a 1-block partition on the circular diagram. The cumulant is assigned to 
the closing edge of the cycle and asymptotic dimensions are assigned to the remaining edges.}
\end{figure}
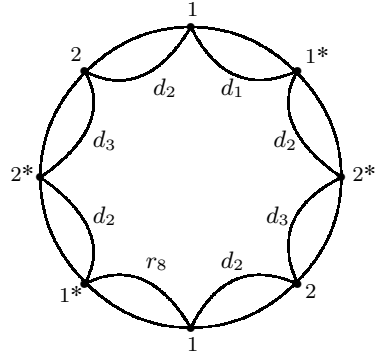

\begin{Example}
{\rm 
For $p=2$ and $k=2$, consider the one-block partition in Fig. 4. It is adapted to the word $W^{2}$, where
$W=122^*1^*$ and all letters have the same label. 
The coloring of edges, discussed in (c) above, is given by $c({\bf 1}_{8})=c_{2}1c_{2}=2321232$, 
where $c_2=232$ and the extended coloring corresponding to the periodic Catalan path is 
$\widetilde{c}({\bf 1}_{8})=123212321$. In the circular diagram we connect all vertices 
of the given block with edges. In this fashion we obtain one additional uncolored edge for each block, to which we assign 
a free cumulant. In this case, it is $r_8$, assigned to the last edge, called the {\it closing edge}
(assuming that the `first' edge is the one connecting 1 and 2 at the bottom of the diagram and we move counterclockwise to reproduce the 
linear diagram). The dimensions are assigned to the remaining edges as shown in Fig. 5. 
}
\end{Example}

\begin{Example}
{\rm For $p=2$ and $k=2$, consider the partition $\pi$ in Fig. 5. It is adapted to the word $W^2$, where
$W=122^*1^*$ and all letters have the same label. It can be represented as the product of three cycles:
$\pi=(1,8)(2,3,4,5)(6,7)$. Note that each of these cycles has an even number of elements and 
each number $s\in [8]$ satisfies either the shift or the reflection condition. For instance, using 
addition modulo $4$, reflection gives $\pi(1)=8=4=4-1+1\;$ and $\pi(8)=1=-3=4-8+1$, 
$\pi(6)=7=-1=4-6+1$ and $\pi(7)=6=-2=4-7+1$, whereas shift gives $\pi(2)=3=2+1, \pi(3)=4=3+1, 
\pi(4)=5=4+1$ and $\pi(5)=2=6=5+1$.}
\end{Example}

\begin{Example}
{\rm 
Consider again the coloring $\widetilde{c}({\bf 1}_{8})=123212321$. 
Let us give three examples of color adapted bracketings and their reductions to the coloring $1$:
$$
1(2(321)2(3)2)1\rightarrow 1(2(3)2)1\rightarrow 1(2)1\rightarrow 1,
$$
$$
1(2(3(212)3)2)1\rightarrow 1(2(3)2)1\rightarrow 1(2)1\rightarrow 1,
$$
$$
1(232)1(2(3)2)1\rightarrow 1(2(3)2)1\rightarrow 1(2)1\rightarrow 1,
$$
where we made one reduction at a time. Each of these bracketings defines a noncrossing partition adapted
to $W^2$ for $p=2$. The partition corresponding to the first bracketing is given by 
$\pi$ in Fig. 4. Altogether, there are 22 adapted colorings in this case and thus 
$|\mathcal{NC}(W^2)|=22$. This number is equal to the 
rescaled Raney number $4^2R_{2}(2,1/2)$ 
and this enumeration involving Raney numbers $R_{k}(2,1/2)$ holds for $p=2$ and all $k$, which will be shown later.}
\end{Example}

\begin{Theorem}
Under assumptions A1-A3 on $\mathpzc{Y}$, let $B=X_1X_2\ldots X_p$, where matrices $X_l=S_{l,l+1}(u_{l})$
are blocks of $Y(u_l)\in \mathpzc{Y}$ for $1\leq l \leq p$, where $u_1, \ldots, u_p\in \mathpzc{U}$.
Then 
$$
\lim_{n\rightarrow \infty}\tau_{1}((BB^*)^{k})=\sum_{\pi\in \mathcal{NC}(W^{k})}
\prod_{{\rm blocks}\, V} w(V),
$$
where the weights assigned to blocks $V$ of $\pi$ are of the form $w(V)=d(V)r_{|V|}(u_V)$.
\end{Theorem}
{\it Proof.}
The proof is similar to that of Corollary 4.1 in [16]. Under all partial traces, and therefore also under 
$\tau_1$, we have convergence
$$
\lim_{n\rightarrow \infty}X_j= \Gamma_{j}\;\;\;{\rm and}\;\;\;
\lim_{n\rightarrow \infty}X_j^*= \Gamma_{j}^*
$$
for any $j$, where 
$\Gamma_{j}=\Gamma_{j,j+1}(u_j)$ and $\Gamma_{j}^{*}=\Gamma_{j+1,j}(u_j).$
Therefore, it suffices to compute the moments 
$$
m_{k}=\Psi_1((\Gamma\Gamma^*)^{k})
$$ 
for $\Gamma=\Gamma_1\ldots \Gamma_p$. Assigning
letters $j,j^*$ to $\Gamma_j, \Gamma_j^*$, respectively, we associate
the word $W^k$ with $m_k$ in the natural fashion. Therefore, 
as in [16], the moment $m_k$ can be expressed in terms of noncrossing partitions $\pi$ of the 
set of $2pk$ numbers which are adapted to the word $W^k$ in the sense specified above. 

The notion of adaptedness to words guarantees that operators corresponding to 
the same block $V$ of $\pi$ have the same label $u$, thus we can write $u=u_{V}$. 
Note that both $\Gamma_j$ and $\Gamma_j^*$ have one annihilation operator, $\wp_{j+1,1}^*$ and $\wp_{j,j+1}^*$, respectively, 
and thus each leg of the block (except the rightmost one) of color $j$ produces color $j+1$ of the edge starting from $j$ and each leg (except the rightmost one)
of color $j^*$ produces the color $j$ of the edge starting from $j^*$, as discussed in (c) above. 
In turn, each block $V$ corresponds to the product of annihilation and creation operators of the type considered in (d) 
above (see Fig. 2), with each annihilation operator taken from a different $\Gamma_{j}$ or $\Gamma_{j}^{*}$ 
and the product of creation operators taken from $\Gamma_{i_{m-1}}$ or $\Gamma_{i_{2}}^*$. 
The contribution of such block is given by the free cumulant 
$r_{k}(u)$ standing by $\wp_{i_{k-1},i_{k-2}}\ldots \wp_{i_{1},j}$ 
multiplied by the product of asymptotic dimensions $d(V)=d_{i_1}d_{i_2}\ldots d_{i_{k-1}}$ obtained from 
the coloring $c(V)=i_1i_2\ldots i_{k-1}$, as discussed in (d). Finally, the property of the whole 
product 
$$
\wp_{i_1,j}^*\wp_{i_2,i_1}^*\ldots \wp_{i_{m-1},i_{m-2}}^{*}\wp_{i_{m-1},i_{m-2}}\ldots \wp_{i_2,i_1}\wp_{i_{1},j}
$$
described in (e) ensures that $c(V)$ is surrounded from both sides by 
color $j$, which entails color adaptedness of $\pi$. This means that each non-vanishing contribution to $m_k$ 
corresponds to some $\pi\in \mathcal{NC}(W^k)$ and it is equal to the product of 
weights $w(V)$ over blocks of $\pi$. Conversely, it is not hard to see that for
any noncrossing partition $\pi$ adapted to $W^k$ we can find the corresponding 
product of operators since legs which are not rightmost determine the places where 
we have to write annihilation operators and their colors are determined by the colorings of blocks, whereas 
the remaining legs must then correspond to products of creation operators and these are uniquely determined by the annihilation 
operators since the partition is noncrossing. This completes the proof.
\hfill $\blacksquare$\\

Two situations are of particular interest and will be studied in more detail in Sections 5-7. 
First, we shall consider the case when $i\neq j\rightarrow u_i\neq u_j$, which corresponds to independent blocks.   
The second case is when $u_j=u$ for any $j\in [p]$, which corresponds to, in general, dependent 
blocks of one random matrix. However, first we will express the limit moments of Theorem 3.1 
in terms of noncrossing pair partitions. This `Gaussianization' holds in the general case.

\section{Gaussianization}
We would like to express the result of Theorem 3.1, using only noncrossing {\it pair} partitions. This will 
produce a  `Gaussianization' of the associated product $BB^*$, by which we mean the fact that
its moments will be expressed in terms of the same class of partitions that is used
when $B$ is a product of Gaussian random matrices and the only difference is that 
we assign some weights to their blocks.

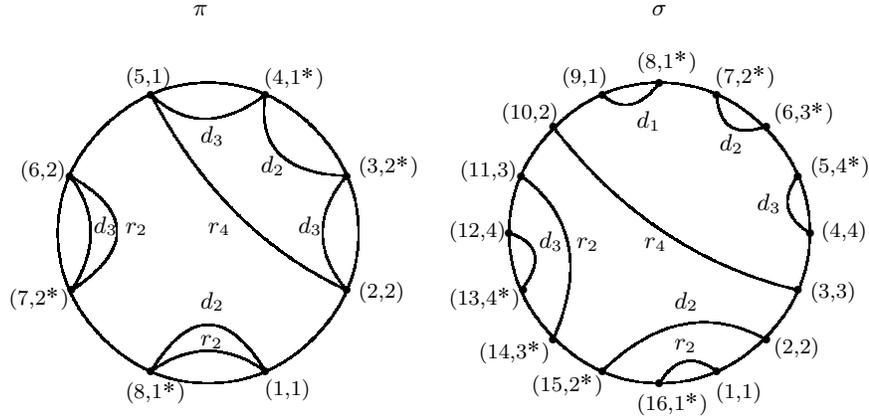
\begin{figure}
\unitlength=1mm
\special{em.linewidth 1pt}
\linethickness{0.5pt}
\begin{picture}(100.00,60.00)(-80.00,-25.00)
%%%%%%%%%%%%%%%%%%%%%%%%%%%%%%%%%%%%%%%%%%%%
\put(-62,29){\scriptsize $\pi$}
\put(-61.00,-15.00){\scriptsize $r_{2}$}
\put(-61.00,-10.00){\scriptsize $d_2$}
\put(-75.20,00.00){\scriptsize $d_3$}
\put(-71.00,00.00){\scriptsize $r_2$}
\put(-60.00,00.00){\scriptsize $r_{4}$}
\put(-61.00,12.00){\scriptsize $d_3$}
\put(-53.00,08.00){\scriptsize $d_2$}
\put(-48.00,00.00){\scriptsize $d_{3}$}
%%%%%% quadrant 1  circle 1 %%%%%%%%%%%%%%%%%%%%%%%%%%%%%
\put(-41.60,07.60){\circle*{1.00}}
\put(-40.00,08.00){$\scriptstyle{(3,2^*)}$}
\put(-52.40,18.40){\circle*{1.00}}
\put(-52.60,19.80){$\scriptstyle{(4,1^*)}$}

%%%%%%%%%%%%%%%%% quadrant 4 circle 1 %%%%%%%%%%%%
\put(-41.60,-07.60){\circle*{1.00}}
\put(-40.00,-8.60){$\scriptstyle{(2,2)}$}
\put(-52.40,-18.40){\circle*{1.00}}
\put(-52.50,-21.60){$\scriptstyle{(1,1)}$}

%%%%%%%%%%%%%%%%% quadrant 3 circle 1 %%%%%%%%%%%%
\put(-67.60,-18.40){\circle*{1.00}}
\put(-71.00,-22.00){$\scriptstyle{(8,1^*)}$}
\put(-78.14,-07.60){\circle*{1.00}}
\put(-86.50,-09.50){$\scriptstyle{(7,2^*)}$}

%%%%%%%%%%%%%%%%% quandrant 2 circle 1 %%%%%%%%%%%
\put(-67.60,18.40){\circle*{1.00}}
\put(-71.00,20.00){$\scriptstyle{(5,1)}$}
\put(-78.40,07.60){\circle*{1.00}}
\put(-85.00,07.50){$\scriptstyle{(6,2)}$}

%%%%%%%%%%%%%%% edges circle 1 %%%%%%%%%%%%%%%%%%%%%%%%%%%%
%%%%%%%%%%%%% 1-1* %%%%%%%%%%%%%%%%%%%%%%%%%%%
\qbezier(-52.40,18.40)(-60.80,12.00)(-67.60,18.40)
%%%%%%%%%%% 2-1 %%%%%%%%%%%%%%%%%%%%%%%%%%%%%%
\qbezier(-41.60,-07.60)(-58.00,00.00)(-67.60,18.40)
%%%%%%%%%%% 2-2* %%%%%%%%%%%%%%%%%%%%%%%%%%%%%%%
\qbezier(-78.14,-07.60)(-73.00,0.50)(-78.40,7.60)
\qbezier(-78.14,-07.60)(-66.00,0.50)(-78.40,7.60)
%%%%%%%%%% 2*-1* %%%%%%%%%%%%%%%%%%%%%%%%%%%
\qbezier(-41.60,07.60)(-54.00,08.00)(-52.40,18.40)
%%%%%%%%%% 1-1* %%%%%%%%%%%%%%%%%%%%%%%%%%%%%%%%%
\qbezier(-67.60,-18.40)(-60.00,-13.00)(-52.40,-18.40)
\qbezier(-67.60,-18.40)(-60.00,-06.00)(-52.40,-18.40)
%%%%%%%%%%% 2- 2* %%%%%%%%%%%%%%%%%%%%%%%%%%%
\qbezier(-41.60,07.60)(-48.00,00.00)(-41.60,-07.60)
%%%%%%%%%%%%%%%%%%%%%%%%%%%%%%%%%%%%%%%%%%%%%%%%%%

%%%%%%%%%%%%%%%%%%% circle 1 from bezier curves %%%%%%%%%%%%%%%%%%%%%%%
\qbezier(-60,20)(-52,20)(-46,14)
\qbezier(-46,14)(-40,08)(-40,00)
\qbezier(-40,00)(-40,-08)(-46,-14)
\qbezier(-46,-14)(-52,-20)(-60,-20)
\qbezier(-60,-20)(-68,-20)(-74,-14)
\qbezier(-74,-14)(-80,-08)(-80,00)
\qbezier(-80,00)(-80,08)(-74,14)
\qbezier(-74,14)(-68,20)(-60,20)

%%%%%%%%%%%%%%%%%%%%%%%%%%%%%%%%%%%%%%%%%%%%
\put(-01,29){\scriptsize $\sigma$}
\put(-02.00,-02.00){\scriptsize $r_{4}$}
\put(-03.00,14.00){\scriptsize $d_{1}$}
\put(02.00,-15.50){\scriptsize $r_{2}$}
\put(13.00,03.00){\scriptsize $d_{3}$}
\put(08.00,11.00){\scriptsize $d_{2}$}
\put(02.00,-10.00){\scriptsize $d_{2}$}
\put(-11.00,-02.00){\scriptsize $r_{2}$}
\put(-16.00,-02.00){\scriptsize $d_{3}$}

%%%%%% quadrant 1  circle 2 %%%%%%%%%%%%%%%%%%%%%%%%%%%%%
\put(00.00,20.00){\circle*{1.00}}
\put(-3.00,22.00){$\scriptstyle{(8,1^*)}$}
\put(20.00,00.00){\circle*{1.00}}
\put(21.50,-00.50){$\scriptstyle{(4,4)}$}
\put(18.40,07.60){\circle*{1.00}}
\put(20.00,08.00){$\scriptstyle{(5,4^*)}$}
\put(14.20,14.20){\circle*{1.00}}
\put(15.20,15.20){$\scriptstyle{(6,3^*)}$}
\put(07.60,18.40){\circle*{1.00}}
\put(07.40,20.00){$\scriptstyle{(7,2^*)}$}

%%%%%%%%%%%%%%%%% quadrant 4 circle 2 %%%%%%%%%%%%
\put(18.40,-07.60){\circle*{1.00}}
\put(20.00,-8.60){$\scriptstyle{(3,3)}$}
\put(14.20,-14.20){\circle*{1.00}}
\put(15.20,-16.00){$\scriptstyle{(2,2)}$}
\put(7.60,-18.40){\circle*{1.00}}
\put(7.50,-22.00){$\scriptstyle{(1,1)}$}
\put(00.00,-20.00){\circle*{1.00}}
\put(-3.00,-23.50){$\scriptstyle{(16,1^*)}$}

%%%%%%%%%%%%%%%%% quandrant 3 circle 2 %%%%%%%%%%%%
\put(-20.00,00.00){\circle*{1.00}}
\put(-28.00,-00.50){$\scriptstyle{(12,4)}$}
\put(-14.14,-14.20){\circle*{1.00}}
\put(-24.00,-16.60){$\scriptstyle{(14,3^*)}$}
\put(-07.60,-18.40){\circle*{1.00}}
\put(-17.00,-21.00){$\scriptstyle{(15,2^*)}$}
\put(-18.14,-07.60){\circle*{1.00}}
\put(-28.00,-09.50){$\scriptstyle{(13,4^*)}$}

%%%%%%%%%%%%%%%%% quandrant 2 circle 2 %%%%%%%%%%%
\put(-07.60,18.40){\circle*{1.00}}
\put(-13.00,20.00){$\scriptstyle{(9,1)}$}
\put(-14.10,14.10){\circle*{1.00}}
\put(-21.00,15.20){$\scriptstyle{(10,2)}$}
\put(-18.40,07.60){\circle*{1.00}}
\put(-26.50,07.50){$\scriptstyle{(11,3)}$}

%%%%%%%%%%%%%%% edges circle 2 %%%%%%%%%%%%%%%%%%%%%%%%%%%%
%%%%%%%%%%%%% 1-1* %%%%%%%%%%%%%%%%%%%%%%%%%%%
\qbezier(00.00,20.00)(-03.80,15.00)(-07.60,18.40)
%%%%%%%%%%% 2-2*%%%%%%%%%%%%%%%%%%%%%%%%%%%%%%
\qbezier(07.60,18.40)(8.00,12.00)(14.20,14.20)
%%%%%%%%%%% 3-3* %%%%%%%%%%%%%%%%%%%%%%%%%%%%%%
\qbezier(-14.14,-14.20)(-08.00,00.00)(-18.40,7.60)
%%%%%%%%%%% 4-4* %%%%%%%%%%%%%%%%%%%%%%%%%%%%%%%
\qbezier(-20.00,0.00)(-14.00,-2.00)(-18.40,-7.60)
%%%%%%%%%%% 3-3* %%%%%%%%%%%%%%%%%%%%%%%%%%%%%%%
\qbezier(-14.20,14.20)(-01.00,-03.00)(18.40,-07.60)
%%%%%%%%%% 2-2* %%%%%%%%%%%%%%%%%%%%%%%%%%%
\qbezier(-07.60,-18.40)(03.00,-08.00)(14.20,-14.20)
%%%%%%%%%% 1-1* %%%%%%%%%%%%%%%%%%%%%%%%%%%%%%%%%
\qbezier(00.00,-20.00)(03.00,-15.00)(07.60,-18.40)
%%%%%%%%%%% 4- 4* %%%%%%%%%%%%%%%%%%%%%%%%%%%
\qbezier(20.00,0.00)(15.00,03.00)(18.40,07.60)
%%%%%%%%%%%%%%%%%%%%%%%%%%%%%%%%%%%%%%%%%%%%%%%%%%

%%%%%%%%%%%%%%%%%%% circle 2 from bezier curves %%%%%%%%%%%%%%%%%%%%%%%
\qbezier(00,20)(08,20)(14,14)
\qbezier(14,14)(20,08)(20,00)
\qbezier(20,00)(20,-08)(14,-14)
\qbezier(14,-14)(08,-20)(00,-20)
\qbezier(00,-20)(-08,-20)(-14,-14)
\qbezier(-14,-14)(-20,-08)(-20,00)
\qbezier(-20,00)(-20,08)(-14,14)
\qbezier(-14,14)(-08,20)(00,20)
%%%%%%%%%%%%%%%%%%%%%%%%%%%%%%%%%%%%%%%%%%%%%%%%%%%%%%%%%%%%%%%%%%%%
\end{picture}
\caption{On the left we have a partition $\pi\in \mathcal{NC}(W^2)$ represented 
in terms of cycles. We assign $(s,j)$ to each vertex, where $s\in [8]$
and $j\in A_2$. 
Dimensions are assigned to regular edges and free cumulants are assigned 
to closing edges, with the convention that the closing edge has odd depth.
On the right we have the corresponding $\sigma=\alpha(\pi)\in \mathcal{NC}^{2}(\widetilde{W}^2)$ represented 
as a pairing $\{\{2s, \beta(s)\}, s\in [8]\}$, where $\beta(s)=2\pi(s)-1$, 
with vertices colored by $A_4$ counterclockwise.
Free cumulants and dimensions are assigned to blocks of odd and 
even depths, respectively.}
\end{figure}

For that purpose, we shall use the bijection between $\mathcal{NC}(n)$ and 
$\mathcal{NC}^{2}(2n)$.  We identify each block $V=\{s_1<s_2<\ldots <s_m\}$ of 
$\pi\in \mathcal{NC}(n)$ with the cycle $(s_1,s_2,\ldots, s_m)$ and set $\pi(s_i)=s_{i+1}$ 
for $1\leq i \leq m-1$ and $\pi(s_m)=s_1$. Then, the mapping
$$
\alpha:\mathcal{NC}(n)\rightarrow \mathcal{NC}^{2}(2n)
$$
given by
$$
\alpha(\pi)=\{\{2s,\beta(s)\}: s\in [n]\},
$$
for $\pi\in \mathcal{NC}(n)$, where 
$$
\beta(s)=2\pi(s)-1
$$
for any $s\in [n]$, gives the required bijection.  

\begin{Remark}
{\rm 
One of the ways to justify that is to compare $\alpha(\pi)$ with the bijection 
$$
\delta(\pi)=\{\{s, \gamma(s)\}:s\in [n]\},
$$
where 
$$
\gamma(s)=2n+1-\pi(s),
$$
used by Mingo and Popa [20]. In the circular diagram for $\delta(\pi)$, 
the numbers from $[n]$ are put in order counterclockwise and they interlace with 
those from $[2n]\setminus [n]$, which are put clockwise, beginning on the right of $n$. 
For our purposes, it is more convenient to compose this mapping with 
the permutation of the set $[2n]$ which puts all numbers $s$ 
in a counterclockwise order on the circle, keeping the diagram intact. This 
leads to $\alpha(\pi)$, which is easy to check on the example shown in Fig. 5. }
\end{Remark}
\begin{Remark}
{\rm 
Equivalently, one can define
$\theta:\mathcal{NC}^{2}(2n)\rightarrow \mathcal{NC}(n)$ as
the restriction of the Kreweras complementation map
$$
K:\mathcal{NC}^{2}(2n)\rightarrow \mathcal{NC}(2n)
$$
to odd numbers, namely 
$$
\theta(\sigma)=K(\sigma)|_{\{1,3,\ldots , 2n-1\}}
$$    
for $\sigma\in \mathcal{NC}^{2}(2n)$. Then $\pi=\theta^{-1}(\sigma)\in \mathcal{NC}(n)$ and thus 
$\alpha$ can be identified  $\theta^{-1}$. For more details on the Kreweras complementation map, see [26].}
\end{Remark}

In order to formulate the main result of this section, we will adapt the bijection $\alpha$ 
to our context of partitions adapted to words. As before, we associate with $BB^*$ the word of length $2p$ of the form
$$
W=12\ldots pp^*\ldots 2^*1^*
$$ 
and then create the twice longer word
$$
\widetilde{W}=12\ldots (2p)(2p)^*\ldots 2^*1^*
$$
whose $m$th power defines the coloring of $n=4pk$ vertices of diagrams corresponding to $\sigma=\alpha(\pi)$.
The consecutive letters in words $W$ and $\widetilde{W}$ 
will be used to color $[n]$ and $[2n]$ counterclockwise, respectively. Both leg coloring maps will be denoted by
$$
s\rightarrow j(s)
$$
where $(s, j(s))\in [n]\times A_{p}$ and $(s,j(s))\in [2n] \times A_{2p}$, respectively.
As concerns labels, we understand that two consecutive colors in $\widetilde{W}$ of the form $2j-1$ and $2j$, 
or $2j^*$ and $2j^*-1$ inherit the same label from color $j$ in $W$.   
Again, the example shown in Fig. 5 serves as an illustration of the coloring maps, with labels again suppressed.

\begin{Example}
{\rm Let us consider $\pi=(1,8)(2,3,4,5)(6,7)\in \mathcal{NC}(8)$ 
of Example 3.2. Let $V_1,V_2,V_3$ be the corresponding cycles (blocks).
Its circular diagram is given by the left diagram in Fig. 5. If we assign 
numbers $s\in [8]$ and colors $j\in A_2$ to the vertices as shown there, 
then $\pi$ is adapted to $W^2$, where $W=122^*1^*$, 
provided all letters have the same label.
Next, we map $\pi$ onto $\sigma=\alpha(\pi)$ by `mapping each vertex 
of $\pi$ onto two neighboring vertices' in the right diagram in Fig. 5 and 
defining the pairing of $\sigma$ by the mapping $s\rightarrow \beta(s)$. 
Pictorially, for instance, the vertex $(1,1)$ is `mapped' onto the pair $\{(1,2),(16,1^*)\}$ 
and the vertex $(2,2)$ is `mapped' onto the pair $\{(3,3),(4,4)\}$.
Computing $\beta(s)$ for all $s\in [8]$ gives then the pair partition $\sigma$. Its subpartitions,
$$
\sigma_1=\{\{1,9\},\{8,16\}\}, \;\;
\sigma_2=\{\{2,14\}, \{3, 13\}, \{4,12\}, \{5,15\}\}, \;\;
\sigma_3=\{\{6,10\}, \{7,11\}\},
$$
can be treated as the images of the blocks $V_1,V_2, V_3$, respectively,  under $\alpha$. 
Finally, we assign colors from the alphabet $A_4$ to the vertices of $\sigma$ as shown in Fig. 5.
Note that this assignment is natural in view of the fact that color $j$ is always `mapped' onto two 
consecutive colors $2j-1$ and $2j$, and color $j^*$ is 
also `mapped' onto two consecutive colors, $(2j)^*$ and $(2j-1)^*$.
We also show how to assign free cumulants and asymptotic dimensions to the edges of $\sigma$.}
\end{Example}

\begin{Remark}
{\rm It should be clear that the image of $\mathcal{NC}(W^k)$ under bijection $\alpha$ depends
on the labels of letters from $A_{p}$.
Denote by $\mathcal{NC}^{2}(\widetilde{W}^{k})$ the set of non-crossing pair 
partitions adapted to the word $\widetilde{W}^{k}$ in the sense that 
each letter $i$ can be paired with $i^*$, but also each letter $2i$ can be paired with 
$2i+1$ and each $(2i+1)^*$ can be paired with $(2i)^*$, provided both letters 
have the same label (recall that in our notation for letters labels are not explicitly shown).
It is clear that $\mathcal{NC}^{2}(\widetilde{W}^k)$ can be identified with $\mathcal{NC}^{2}(W_{0}^{k})$, 
where 
$$
W_{0}=122^*44^*\ldots 2p(2p)^*\ldots 44^*22^*1^*,
$$ 
in which all pairings are colored by pairs of the form $\{j,j^*\}$, where $j=1,2,4,\ldots , 2p$, 
provided $j$ and $j^*$ inherited the same label. This is a priori not guaranteed since in this 
realization certain legs colored by $j$ and $j^*$ may inherit different labels. In that respect, using $W_{0}$
with suppressed labels is slightly more delicate. 
}
\end{Remark}
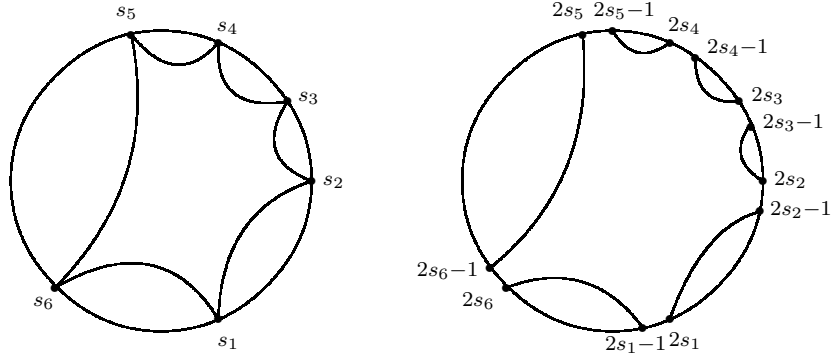
\begin{figure}
\unitlength=1mm
\special{em.linewidth 1pt}
\linethickness{0.5pt}
\begin{picture}(120.00,60.00)(-90.00,-25.00)
%%%%%%%%%%%%%%%%%%%%%%%%%%%%%%%%%%%%%%%%%%%%
%%%%%%%%%%%%%%%%%%%%%%%%%%%%%%%%%%%%%%%%%%%%
%%%%%% vertices %%%%%%%%%%%%%%%%%%%%%%%%%%%%%
\put(-64.00,19.50){\circle*{1.00}}
\put(-66.00,22.00){$\scriptstyle{s_5}$}
\put(-52.40,-18.40){\circle*{1.00}}
\put(-52.50,-22.00){$\scriptstyle{s_1}$}
\put(-40.00,00.00){\circle*{1.00}}
\put(-38.50,-00.50){$\scriptstyle{s_2}$}
\put(-43.20,10.60){\circle*{1.00}}
\put(-41.80,10.60){$\scriptstyle{s_3}$}
\put(-52.40,18.40){\circle*{1.00}}
\put(-52.60,20.00){$\scriptstyle{s_4}$}
\put(-74.14,-14.20){\circle*{1.00}}
\put(-77.00,-16.60){$\scriptstyle{s_6}$}
%%%%%%%%%%%%%%%%%% edges %%%%%%%%%%%%%%%%%%%%%%%%%%%%%
\qbezier(-74.14,-14.20)(-60.00,0.00)(-64.00,19.50)
\qbezier(-52.40,18.40)(-58.00,12.00)(-64.00,19.50)
\qbezier(-40.00,0.00)(-48.00,03.00)(-43.20,10.60)
\qbezier(-52.40,18.40)(-53.00,9.00)(-43.20,10.60)
\qbezier(-52.40,-18.40)(-52.00,-4.00)(-40.00,00.00)
\qbezier(-52.40,-18.40)(-60.00,-6.00)(-74.14,-14.20)
%%%%%%%%%%%%%%%%%%% circle from bezier curves %%%%%%%%%%%%%%%%%%%%%%%
\qbezier(-60,20)(-52,20)(-46,14)
\qbezier(-46,14)(-40,08)(-40,00)
\qbezier(-40,00)(-40,-08)(-46,-14)
\qbezier(-46,-14)(-52,-20)(-60,-20)
\qbezier(-60,-20)(-68,-20)(-74,-14)
\qbezier(-74,-14)(-80,-08)(-80,00)
\qbezier(-80,00)(-80,08)(-74,14)
\qbezier(-74,14)(-68,20)(-60,20)
%%%%%%%%%%%%%%%%%%%%%%%%%%%%%%%%%%%%%%%%%%%%%%%%%%%%%%%%%%%%%%%%%%%%
%%%%%% vertices %%%%%%%%%%%%%%%%%%%%%%%%%%%%%
\put(00.00,20.00){\circle*{1.00}}
\put(-2.50,22.00){$\scriptstyle{2s_5-1}$}
\put(-04.00,19.50){\circle*{1.00}}
\put(-8.00,22.00){$\scriptstyle{2s_5}$}
\put(7.60,-18.40){\circle*{1.00}}
\put(7.50,-22.00){$\scriptstyle{2s_1}$}
\put(20.00,00.00){\circle*{1.00}}
\put(21.50,-00.50){$\scriptstyle{2s_2}$}
\put(19.60,-04.00){\circle*{1.00}}
\put(21.00,-04.50){$\scriptstyle{2s_2-1}$}
\put(16.80,10.60){\circle*{1.00}}
\put(18.50,10.60){$\scriptstyle{2s_3}$}
\put(18.40,07.20){\circle*{1.00}}
\put(20.00,07.00){$\scriptstyle{2s_3-1}$}
\put(07.60,18.40){\circle*{1.00}}
\put(07.40,20.00){$\scriptstyle{2s_4}$}
\put(11.00,16.40){\circle*{1.00}}
\put(12.60,17.00){$\scriptstyle{2s_4-1}$}
\put(-14.14,-14.20){\circle*{1.00}}
\put(-20.00,-16.60){$\scriptstyle{2s_6}$}
\put(-16.30,-11.60){\circle*{1.00}}
\put(-26.00,-12.60){$\scriptstyle{2s_6-1}$}
\put(04.00,-19.50){\circle*{1.00}}
\put(-01.00,-22.50){$\scriptstyle{2s_1-1}$}
%%%%%%%%%%%%%%%%%% edges %%%%%%%%%%%%%%%%%%%%%%%%%%%%%
\qbezier(-16.30,-11.60)(-02.00,0.00)(-04.00,19.50)
\qbezier(07.60,18.40)(03.00,15.00)(00.00,20.00)
\qbezier(-14.14,-14.20)(-03.00,-10.00)(04.00,-19.50)
\qbezier(7.60,-18.40)(12.00,-06.00)(19.80,-04.00)
\qbezier(20.00,0.00)(15.00,03.00)(18.40,07.60)
\qbezier(11.00,16.40)(11.00,10.00)(16.80,10.60)
%%%%%%%%%%%%%%%%%%% circle from bezier curves %%%%%%%%%%%%%%%%%%%%%%%
\qbezier(00,20)(08,20)(14,14)
\qbezier(14,14)(20,08)(20,00)
\qbezier(20,00)(20,-08)(14,-14)
\qbezier(14,-14)(08,-20)(00,-20)
\qbezier(00,-20)(-08,-20)(-14,-14)
\qbezier(-14,-14)(-20,-08)(-20,00)
\qbezier(-20,00)(-20,08)(-14,14)
\qbezier(-14,14)(-08,20)(00,20)
%%%%%%%%%%%%%%%%%%%%%%%%%%%%%%%%%%%%%%%%%%%%%%%%%%%%%%%%%%%%%%%%%%%%
\end{picture}
\caption{Cycle $V=(s_1,s_2,s_3,s_4,s_5,s_6)$ of $\pi$ and the corresponding pair 
subpartition $\sigma'=\{(2s_i,2s_{i+1}-1),i\in [6]\}$ of $\sigma=\alpha(\pi)$. Observe that 
$\sigma'$ can be viewed as a continuous deformation of $V$ which maps each vertex of $V$ onto 
a pair of neighboring vertices of $\sigma$.}
\end{figure}

\begin{Lemma}
There is a bijection between $\mathcal{NC}(W^{k})$ and $\mathcal{NC}^{2}(\widetilde{W}^{k})$.
\end{Lemma}
{\it Proof.}
By definition, 
$$
\mathcal{NC}(W^{k})\subset \mathcal{NC}(n)\;\;\;{\rm and} \;\;\;
\mathcal{NC}^2(\widetilde{W}^k)\subset \mathcal{NC}^{2}(2n)
$$ 
where $n=2pk$. It is well known that $\mathcal{NC}(n)\cong \mathcal{NC}^{2}(2n)$, therefore
it suffices to show two implications: $\alpha(\pi)\in \mathcal{NC}^{2}(\widetilde{W}^{k})$ 
for any $\pi\in \mathcal{NC}(W^{k})$,
and $\alpha^{-1}(\sigma)\in \mathcal{NC}(W^k)$ for any $\sigma\in \mathcal{NC}^{2}(\widetilde{W}^{k})$.\\
\indent
{\it 1st implication.}
If $\pi\in \mathcal{NC}(W^{k})$, we know that $\alpha(\pi)$ 
consists of pairs of numbers of the form 
$$
\alpha(\pi)=\{\{2s, 2\pi(s)-1\}: s\in [p]\},
$$
of which the first one is even and the second one (smaller or bigger than $2s$) 
is odd.
Let us restrict our attention to the case when $s$ belongs 
to the block of $\pi$ of the form $V=(s_1,s_2, \ldots ,s_m)$. 
Recall that by an edge of $V$ we understand a pair 
$(s_i,s_{i+1})$, where $i\in [m]$ and $s_{m+1}=s_{1}$.
For instance, the first diagram in Fig. 6 has six edges.
Further, since $\pi$ is adapted to $W^{k}$, there are only four types of colorings of vertices 
corresponding to edges:
$$
E=\{(j,j+1),(j,j^*),(j^*,j^{*}-1),(j^*,j)\}.
$$ 
Of course, since vertices are colored by $A_{p}$ in a natural counterclockwise 
order, each of these four colorings is possible only for selected edges. 
For instance, if the pair $(s,s')$ is colored by $(j,j+1)$, then
$s'=s+1\;({\rm mod}\;2p)$, if $(s,s')$ is colored by $(j,j^*)$, then $s'=2p-s\;({\rm mod}\; 2p)$, etc.
These relations can be easily determined for all four types of colorings, but it suffices to
know that there are four admissible colorings of each pair $(s,s')$ in order to carry out the proof. 
Namely, the mapping $\alpha$ induces
the mapping on the set of colored vertices $(s,j)$ of $\pi$ of the form
$$
\alpha((s,j))=(2s,2j)\;\;\;{\rm and}\;\;\;\alpha((s,j^*))=(2s,2j^*-1),
$$
where we understand that $2j^*-1=(2j-1)^*$ and, by abuse of notation, we also use the symbol $\alpha$ 
for the induced mapping. To complete the proof that $\sigma=\alpha(\pi)$ is adapted to $\widetilde{W}^{k}$,
we need to analyze all edges of $\sigma$ and the colorings of their vertices. 
Let $(s,s')$ be an edge in some $V$. One can think of $(2s, 2s'-1)$ as the corresponding edge in $\sigma'$.
Again, we refer the reader to Fig. 6, where there are six such edges in the right diagram.
We would like to determine the colorings of the vertices of $\sigma'$ which form edges.
There are 4 possible cases of coloring these vertices, depending on the coloring $(j,j')$ of $(s,s')$, namely
\begin{eqnarray*}
(j,j+1)&\rightarrow& (2j,2j+1)\\
(j,j^*)&\rightarrow& (2j,2j^*)\\
(j^*,j)&\rightarrow& (2j^*-1,2j-1)\\
(j^*,j^*-1)&\rightarrow& (2j^*-1,2j^*-2).
\end{eqnarray*}
For instance, if $(s,s')$ is an edge of $V$ and $s$ and $s'$ are colored by $j$ and $j+1$, respectively, 
then we treat $(2s, 2s'-1)$ as the corresponding edge in $\sigma'$. The coloring of $2s$ is $2j$ and 
the coloring of $2s'-1$ is $2j'-1=2j+1$ since $j'=j+1$, which gives the coloring of 
the vertices of $\sigma'$ in the first case. The remaining three cases are treated in a similar way. 
This shows that the coloring of $\sigma'$ is adapted to $\widetilde{W}^{k}$ and thus $\sigma=\alpha(\pi)$ is 
adapted to $\widetilde{W}^{k}$.\\
\indent
{\it 2nd implication.}
In this part of the proof we will use the notion of a noncrossing partition of even depth.
We will say that a noncrossing partition $\sigma$ has {\it even depth} if each deepest 
block of $\sigma$ (block that does not have any inner blocks) has even depth. 
For instance, the partition of Fig. 1 has even depth. 
We will show first that each $\sigma\in \mathcal{NC}^2(\widetilde{W}^{k})$ is of even depth.
Note that the deepest blocks of $\sigma$ must connect neighboring numbers and thus they must be colored 
by $(2j,2j+1)$, $(2p, 2p^*)$ or $(2j^*-1,2j^*-2)$ for some $j\in \{1, \ldots , p-1\}$ 
and $j^*\in \{2^*, \ldots p^*\}$. In the case of $(2j, 2j+1)$, there is an odd number of legs on the left of $2j$ 
and on the right of $2j+1$. If there are some blocks with both left and right legs lying on 
the same side of the pair $(2j,2j+1)$, then the remaining legs lying on that side are either all left 
legs or all right legs and there is still an odd number of them. 
Therefore, the depth of $(2j,2j+1)$ must be even. A similar argument works for the other two cases, which proves that 
$\sigma$ is of even depth. This allows us to decompose $\sigma$ into a collection of 
subpartitions $\sigma'$ consisting of a block of odd depth and all its nearest inner blocks.
It is clear that all left (right) legs in each $\sigma'$ must be associated with odd (even) numbers $s\in [2n]$, where $n=2pk$.
Therefore, $\sigma'$ is of the form 
$$
\sigma'=\{(2s_i, 2s_{i+1}-1): i\in [m]\}
$$
where we set $s_{m+1}=s_{1}$ (as in the right diagram in Fig. 6). It is not difficult to see that $\sigma'$ must have
an even number of blocks. For, if the left leg of the outer block is colored by $j\in [p-1]$ 
(the right leg is then colored by $j^*$), its nearest inner blocks must be colored by alternating 
pairs $(j+1, j^*+1)$ and $(j^*,j)$, the last one being colored again by $(j+1, j^*+1)$, thus the  number of nearest inner blocks
must be odd, so the total number of blocks of such $\sigma'$ must be even. An analogous reasoning holds if 
the left leg of the outer block of $\sigma'$ is colored by $j^*+1$, with the right leg colored by $j+1$.
On this occasion observe that blocks colored by the pair $(2p,2p^*)$ must have even depth and need not be considered 
at this moment. Therefore, $\sigma'$ is exactly of the same form as pair partitions which are 
the images of blocks $V$ of $\pi\in \mathcal{NC}(W^{k})$ under $\alpha$. The colors assigned to 
the legs of edges of $\sigma'$ are 
of type $(2j,2j-1)$, $(2j^*-1,2j^*-2)$, $(2j,2j^*)$ or $(2j^*,2j)$. Therefore, 
the partition $\pi=\alpha^{-1}(\sigma)$, considered as an element of 
$\mathcal{NC}(n)$, consists of blocks $V=(s_{1}, \ldots , s_{k})$ such that 
each pair of neighboring legs $(s_i,s_{i+1})$ is colored by $(j,j+1)$, $(j^*-1,j^*-2)$, $(j,j^*)$ or $(j^*,j)$, respectively. However, this means that $\pi$ is adapted to $W^{k}$ since it acts on each leg $s_i$ 
as a shift or a reflection (modulo $2p$). This completes the proof of the second implication.
\hfill $\blacksquare$\\

Using Lemma 4.1, we would like to express Theorem 3.1 in terms of noncrossing pair partitions
$\sigma \in \mathcal{NC}^{2}(\widetilde{W}^{m})$.
For that purpose, we need to determine the right dimension function for blocks of $\sigma$.
It will suffice to assign dimensions only to blocks of even depths. It is convenient to do it, using 
the alphabet for $W_{0}$. In that setting, the legs of all blocks are colored by natural pairs: 
$(1,1^*)$, $(1^*,1)$, $(2j,2j^*)$ or $(2j^*,2j)$.
We do it as follows:\\
\unitlength=1mm
\special{em:linewidth 1pt}
\linethickness{0.4pt}
\begin{picture}(85.00,20.00)(-60.00,-4.00)
%%%%%%%%%%%%%%%%%%%%%%%%%%%%%%%%%%%%%%%%%%%%%%%%%%%%%%%
\put(-32.50,5.00){$d($}
\put(-28.00,4.00){\line(0,1){4.00}}
\put(-22.00,4.00){\line(0,1){4.00}}
\put(-28.00,8.00){\line(1,0){6.00}}
\put(-21.00,5.00){$)=d_{1},$}
\put(-29.00,1.00){$\scriptscriptstyle{1^*}$}
\put(-23.00,1.00){$\scriptscriptstyle{1}$}

\put(42.00,4.00){\line(0,1){4.00}}
\put(48.00,4.00){\line(0,1){4.00}}
\put(42.00,8.00){\line(1,0){6.00}}

\put(-2.50,5.00){$d($}
\put(2.00,4.00){\line(0,1){4.00}}
\put(8.00,4.00){\line(0,1){4.00}}
\put(2.00,8.00){\line(1,0){6.00}}
\put(9.00,5.00){$)=d_{j+1}$}
\put(01.00,1.00){$\scriptscriptstyle{2j^*}$}
\put(7.00,1.00){$\scriptscriptstyle{2j}$}

\put(42.00,4.00){\line(0,1){4.00}}
\put(48.00,4.00){\line(0,1){4.00}}
\put(42.00,8.00){\line(1,0){6.00}}

\put(27.00,5.00){{\rm and}}
\put(37.50,5.00){$d($}
\put(41.00,1.00){$\scriptscriptstyle{2j}$}
\put(47.00,1.00){$\scriptscriptstyle{2j^*}$}
\put(49.00,5.00){$)=d_{j+1}$}
\end{picture}
\\
whenever these blocks are of even depths. Let us note that each block colored by $(1,1^*)$ is always of 
odd depth and that is why we did not assign to it any dimension.
\begin{Theorem}
Under the assumptions of Theorem 3.1, it holds that
$$
\lim_{n\rightarrow \infty}\tau_{1}((BB^*)^{k})=\sum_{\sigma\in \mathcal{NC}^2(W^{k}_{0})}
\prod_{{\rm blocks}\, V\; {\rm of}\;\sigma} w(V),
$$
where the weights assigned to blocks $V$ of $\sigma$ are of the form 
$$
w(V)=\left\{
\begin{array}{ll}
r_{\mathpzc{i}(V)+1}(u_V)& if \;\;\mathpzc{d}(V)\;odd\\
d(V)&if \;\;\mathpzc{d}(V)\; even
\end{array}
\right.
$$
and $d(V)$ is assigned to blocks as above.
\end{Theorem}
{\it Proof.}
We use Theorem 3.1 and Lemma 4.1. 
Let us observe that each block $V_1$ of $\pi\in \mathcal{NC}(W^{k})$ consisting of $2k$ elements 
is mapped by $\alpha$ onto a subpartition $\sigma'$ of $\sigma$ of depth two consisting of one outer block 
and its $2k-1$ nearest inner blocks.  Since $r_{2k}(u)$ is assigned to $V_1$, we choose to assign
this cumulant to the distinguished outer block $V$ of $\sigma'$, and therefore $r_{2k}(u)=r_{\mathpzc{i}(V)+1}(u)$.
It remains to check how to assign dimensions to its nearest inner blocks.
Thus, block of $\sigma'$ colored by $(1^*,1)$ corresponds to a segment of $\pi'=\alpha^{-1}(\sigma')$ with legs 
colored by $(1^*,1)$ and the dimension
corresponding to such segment was $d_1$. Similarly, block colored by $(2j,2j^*)$ corresponds to a segment of
$\pi'$ with legs colored by $(j,j+1)$ or $(j,j^*)$ and the corresponding dimension was in both cases $d_{j+1}$. 
Finally, block colored by $(2j^*,2j)$ corresponds to a segment of $\pi'$ with legs colored by $(j^*+1,j^*)$
and the corresponding dimension was $d_{j+1}$. The product of all these dimension weights 
is the same as that assigned to $V'$. Therefore, 
$d(V_1)=\prod_{{\rm blocks}\;V\;{\rm of}\; \sigma'}d(V)$, which completes the proof.
\hfill $\blacksquare$.\\

We view the above theorem as a `Gaussianization' of the asymptotic distribution. 
For each $m$, we have the same class of partitions as in the case when $B$ is a product of $2p$ 
independent GRM and we only need to assign certain weights to the blocks of these partitions.
In particular, if we set all dimensions to be equal to one, the weight assigned to 
a block depends only on the number of its nearest inner blocks. 

\begin{Example}
{\rm Consider again $\pi\in \mathcal{NC}(W^{2})\subset \mathcal{NC}(8)$  
given in Fig. 5, with all letters having the same label. It consists of 3 blocks: $V_1=(1,8)$, $V_2=(2,3,4,5)$, $V_3=(6,7)$.
The corresponding weights are $w(V_1)=r_2d_2$, $w(V_2)=r_4d_1d_2d_3$, $w(V_3)=r_2d_3$, 
and thus the contribution of $\pi$ is 
$$
w(\pi)=r_{4}r_{2}^{2}d_{1}d_{2}^{2}d_{3}^{2}.
$$ 
In particular, if all asymptotic dimensions are equal to $d$, it is of order $d^{5}$.
This term does not show up in the standard Wishart case, when each $Y(u_j)$ is Gaussian, since in 
that case free cumulants of orders higher than two vanish (in the Gaussian case, all
terms which give contribution to $M_2$ are of order $d^{4}$). The corresponding pair partition 
$\alpha(\pi)\in \mathcal{NC}^{2}(\widetilde{W}^{2})\subset \mathcal{NC}^{2}(16)$ 
shown in Fig. 5 consists of 8 pairs. Note that $V_2$ is mapped by $\alpha$ onto 2 pairs and 
the remaining blocks are pairs, thus each of them is mapped by $\alpha$ onto 2 pairs. We place
the free cumulants next to arcs of odd depths, between this arc and the arcs which correspond to
the inner blocks (in order to find the depth of a given block on a circular diagram,
one has to choose the first number, and we understand that it is the 
bottom $1$).}
\end{Example}

\begin{Example}
{\rm In contrast to the previous example, take $p=2$ and $m=2$ and assume that $1$ and $2$ 
have different labels (they correspond to blocks of independent matrices $Y_1, Y_2$).
Then is is not hard to see that the set $\mathcal{NC}(W^{2})$ consists of 
$5$ partitions. Theorem 3.1 gives 
$$
M_{2}=r_{4}s_{2}^2d_1d_2^2d_3^2+r_{2}^2s_{4}d_1d_2^2d_3^2+
r_{2}^2s_{2}^2 (d_1d_2d_3^2+d_2^2d_3^2++d_1d_2^2d_3)
$$
where $r_4,r_2$ and $s_4,s_2$ are free cumulants associated with $Y_1$ and $Y_2$, respectively.
If all dimensions are equal to one, this simplifies to 
$$
M_{2}=r_{4}s_2^2+r_{2}^2s_4+3r_{2}^2s_{2}^2 
$$
and if all free cumulants are equal to one (thus both matrices, $Y(u_1)$ and $Y(u_2)$, from which we take 
our blocks, are standard Wishart matrices), then $M_2=5=F_{2}(4)=|\mathcal{NC}(\widetilde{W}^{2})|$, where the last equation 
follows from Lemma 4.1.}
\end{Example}

\section{Products of independent blocks}

Let us study the generating functions associated with the limit moments 
$$
m_k = \lim_{n\rightarrow\infty} \tau_1((BB^*)^k),
$$
in the case when $X_1, \ldots , X_p$ are consecutive blocks lying above the main diagonal and 
taken from independent random matrices $Y_1, \ldots , Y_p$, where $Y_j=Y(u_j)$ for pairwise 
different $u_j$.
The matrices $X_jX_j^*$ are natural generalizations of the classical Wishart matrix, obtained when 
matrices $Y_j$ are Gaussian.
The associated probability distribution 
on the nonnegative real axis will be expressed in terms of 
even free cumulants of $\nu_1, \ldots , \nu_p$, where $\nu_j$
is the asymptotic distribution of $Y_j$ under $\tau$, using Theorem 3.1. This will lead 
to a formula for moments $m_k$, expressed as linear combinations of generalized 
multivariate Fuss-Narayana polynomials.

For distribution $\mu$ with moments $(m_{k})$ we denote by $S_{\mu}$ the associated $S$-transform 
defined by the formula
$$
S_{\mu}(z)=\frac{z+1}{z}\psi^{-1}_{\mu}(z)
$$
where $\psi_{\mu}$ is the moment generating function without constant term of the form
$$
\psi_{\mu}(z)=\sum_{k=1}^{\infty}m_{k}z^{k}
$$
and $\psi^{-1}$ is its composition inverse. We will also use the reciprocal 
$S$-transform, namely 
$$
T_{\mu}(z)=\frac{1}{S_{\mu}(z)},
$$ 
called the $T$-{\it transform} of $\mu$, used by Dykema [9] and Nica [25]. 
We will use these 
objects also when they are only formal power series associated with 
noncommutative random variables with distribution $\mu$.

We will also use two transformations of distributions, $U_{s}$ and $V_{s}$, for positive $s$. 
For distribution $\mu$ and positive $s$, let $U_{s}$ be defined by
$$
G_{U_{s}(\mu)}(z)=sG_{\mu}(z)+(1-s)z^{-1}
$$
and $G_{\mu}$ is the Cauchy transform of $\mu$. Another useful 
transformation will be defined in terms of $S$-transforms by
$$
S_{V_{s}(\mu)}(z)=S_{\mu}(s^{-1}z)
$$
for any positive $s$. It is easy to see that 
\begin{eqnarray*}
U_{s}(\mu\boxtimes \nu)&=&U_{s}(\mu)\boxtimes V_{s}(\nu)\\
V_{s}(\mu\boxtimes \nu)&=&V_{s}(\mu)\boxtimes V_{s}(\nu)
\end{eqnarray*}
for any distributions $\mu, \nu$ and positive $s$, where $\boxtimes$ is the multiplicative free convolution. 
Moreover, it is clear that
$V_{s}V_{t}=V_{st}$ for any positive $s,t$.

We begin with a generalization of Theorem 9.2 in [16], where rectangular Gaussian random 
matrices were treated and $\mu_1$ was the multiplicative convolution of 
dilated Marchenko-Pastur distributions. In the present context, we replace the latter 
by more general distributions.  
\begin{Proposition}
Under the assumptions of Theorem 3.1, if $X_1, \ldots, X_p$ are blocks of 
independent random matrices $Y(u_1), \ldots , Y(u_p)$, respectively, then
$$
V_{d_1}(\mu)=V_{d_{1}}(\xi_{1})\boxtimes V_{d_{2}}(\xi_{2})\boxtimes \ldots \boxtimes V_{d_{p}}(\xi_{p}),
$$
where $\mu$ is the asymptotic distribution of $BB^*$ under $\tau_1$ 
and $\xi_i$ is the asymptotic distribution of $X_{i}X_{i}^{*}$ under $\tau_i$ for $i\in [p]$.
\end{Proposition}
{\it Proof.}
The proof is similar to that of Theorem 9.2 in [16]. 
Thus, if $\theta_i$ and $\mu_i$ are the asymptotic distributions of $X_{i}^*X_{i}$ and 
$X_{i+1}\ldots X_{p}X_{p}^*\ldots X_{i+1}^*$ under $\tau_{i+1}$, respectively, then
\begin{eqnarray*}
\mu_{1}&=&U_{s_{1}}(\theta_{1}\boxtimes \mu_2),\\
\mu_{2}&=&U_{s_{2}}(\theta_{2}\boxtimes \mu_3),\\
&\ldots&\\
\mu_{p-1}&=&U_{s_{p-1}}(\theta_{p-1}\boxtimes \mu_{p}),
\end{eqnarray*}
where $\mu=\mu_1$ and $s_{i}=d_{i+1}/d_{i}$. 
Now, observe that $\xi_{i}=U_{s_{i}}(\theta_{i})$ for $i\in [p-1]$ and 
$\xi_{p}=\mu_{p}$. In order to use these distributions, it suffices to use 
the properties of transformations $U_s,V_s$, which gives
\begin{eqnarray*}
\mu_{1}&=&\xi_1\boxtimes V_{s_1}(\mu_2),\\
\mu_{2}&=&\xi_2\boxtimes V_{s_2}(\mu_3),\\
&\ldots&\\
\mu_{p-1}&=&\xi_{p-1}\boxtimes V_{s_{p-1}}(\mu_{p}).
\end{eqnarray*}
Since $s_{1}\ldots s_{k}=d_{k+1}/d_{1}$, this yields
$$
\mu=\xi_{1}\boxtimes V_{d_{2}/d_1}(\xi_{2})\boxtimes \ldots \boxtimes V_{d_{p}/d_1}(\xi_{p}),
$$
which is equivalent to the (more symmetric) formula that we needed to prove. 
\hfill $\blacksquare$\\

\begin{Remark}
{\rm Let us observe that in the case of Gaussian random matrices, the asymptotic distributions 
$\xi_{i}$ become dilations of Marchenko-Pastur laws. 
Using the notation of [16, Theorem 9.2], we have
$\xi_{i}=\rho_{d_{i+1},d_{i}}$, with $S$-transforms of the form
$$
S_{\xi_{i}}(z)=(d_{i+1}+d_{i}z)^{-1}
$$
where $i\in [p]$. Thus $\xi_i$ is the $d_{i}$-dilation of the Marchenko-Pastur law 
$\rho_{d_{i+1}}$, where $S_{\rho_{d}}(z)=(d+z)^{-1}$.  Writing Proposition 5.1 in 
terms of $S$-transforms, we obtain 
$$
S_{\mu}(z)=(d_{2}+d_{1}z)^{-1}(d_{3}+d_{1}z)^{-1}\ldots (d_{p+1}+d_{1}z)^{-1}
$$
the formula given in [16, Theorem 9.2]. Note that the computation used in the proof of that theorem 
was slightly different since we kept $\theta_1, \ldots , \theta_{p-1}$, where $\theta_i=\rho_{d_i, d_{i+1}}$. 
In the general case, it is more convenient to use $\xi_1, \ldots , \xi_{p}$ 
since this leads to 
$$
S_{\mu}(d_{1}^{-1}z)
=
S_{\rho_1}(d_1^{-1}z)
\ldots 
S_{\rho_p}(d_p^{-1}z),
$$
a symmetric formula for the $S$-transform of $\mu$.
}
\end{Remark}

The last proposition gives us a relation between $\mu$ and distributions $\xi_1, \ldots , \xi_{p}$.
We would like to find a relation between $\mu$ and distributions 
$\nu_1, \ldots, \nu_p$. The simplest Gaussian counterpart of this relation is 
obtained when $\mu$ is the Marchenko-Pastur law as the limit spectral distribution of $XX^*$, where
$X$ is an off-diagonal block of $Y$, and $\nu$ is the Wigner semicircle law 
as the limit spectral distribution of $Y$. 

In order to do that, we shall use only even free cumulants of $\nu_1, \ldots , \nu_p$.
In general, if $\nu$ is a distribution with free cumulants $(r_{n})$, we denote 
by ${\widetilde{\nu}}$ the distribution given by 
the formal power series (combinatorial $R$-series) of the form
$$
R_{\widetilde{\nu}}(z)=\sum_{n=1}^{\infty}r_{2n}z_{n},
$$
thus $\widetilde{\nu}$ is the distribution whose free cumulants are even free cumulants of the 
probability distribution $\nu$.

\begin{Remark}
{\rm In general, $\widetilde{\nu}$ is not a probability distribution. For instance, if we take $\nu$ to be the free Meixner law 
with parameters $u=v$ and $0<b<a$, then 
$$
r_{2n} = \frac{a(b - a)^{n - 1}}{n}\binom{2n - 2}{n - 1}
$$
and its odd free cumulants vanish. It is easy to check that $(r_{2n})_{n=1}^\infty$ cannot be a free cumulant sequence 
of any probability measure, since the corresponding moment sequence is not positive definite. 
It seems to be an interesting open problem to find necessary and sufficient conditions under which $\widetilde{\nu}$ 
defines a probability measure (and conditions under which it is unique). 
Nevertheless, one can show that if $\nu$ is freely infinitely divisible, then $\widetilde{\nu}$ defines a probability measure. 
This observation as well as the above counterexample was communicated to us by W. M{\l}otkowski.}
\end{Remark}

\begin{figure}
\unitlength=1mm
\special{em:linewidth 1pt}
\linethickness{0.5pt}
\begin{picture}(120.00,30.00)(-05.00,0.00)

\put(50.00,21.00){$\scriptstyle{V_0}$}

\put(10.00,10.00){\line(0,1){8}}
\put(10.00,18.00){\line(1,0){84}}
\put(10.00,10.00){\circle*{1.00}}
\put(17.00,11.00){$\scriptstyle{\pi_1}$}
\put(9.00,05.00){$\scriptstyle{1}$}
\put(16.00,5.00){$\scriptstyle{W_*^{n_1}}$}

\put(28.00,10.00){\line(0,1){8}}
\put(28.00,10.00){\circle*{1.00}}
\put(35.00,11.00){$\scriptstyle{\pi_2}$}
\put(27.00,05.00){$\scriptstyle{1^*}$}
\put(34.00,5.00){$\scriptstyle{W^{n_2}}$}

\put(46.00,10.00){\line(0,1){8}}
\put(58.00,10.00){\line(0,1){8}}
\put(46.00,10.00){\circle*{1.00}}
\put(58.00,10.00){\circle*{1.00}}
\put(63.50,11.00){$\scriptstyle{\pi_{2m-2}}$}
\put(45.00,05.00){$\scriptstyle{1}$}
\put(57.00,05.00){$\scriptstyle{1^*}$}
\put(51.00,5.00){$\scriptstyle{\ldots}$}
\put(62.50,5.00){$\scriptstyle{W^{n_{2m-2}}}$}

\put(76.00,10.00){\line(0,1){8}}
\put(76.00,10.00){\circle*{1.00}}
\put(81.50,11.00){$\scriptstyle{\pi_{2m-1}}$}
\put(75.00,05.00){$\scriptstyle{1}$}
\put(80.50,5.00){$\scriptstyle{W_*^{n_{2m-1}}}$}

\put(94,10.00){\line(0,1){8}}
\put(94,10.00){\circle*{1.00}}
\put(102.50,11.00){$\scriptstyle{\pi_{2m}}$}
\put(93.00,05.00){$\scriptstyle{1^*}$}
\put(101.50,5.00){$\scriptstyle{W^{n_{2m}}}$}

\end{picture}
\caption{Decomposition of $\pi\in \mathcal{NC}^{2}(W^{k})$}
\end{figure}
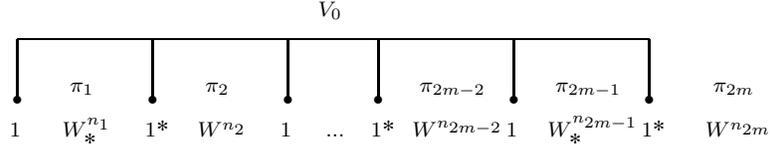

Let us show now that the deformation of the asymptotic distribution of $X_{i}X_{i}^*$ is 
obtained from the corresponding asymptotic distribution of $Y_i$ by convolving 
it with a Marchenko-Pastur distribution.

\begin{Lemma}\label{l1}
Under the assumptions of Proposition 5.1, it holds that 
$$
V_{d_i}(\xi_i)=\widetilde{\nu}_{i}\boxtimes \rho_{d_{i+1}}
$$
for $i\in [p]$, where $\widetilde{\nu}_{i}$ is 
defined by the even free cumulants of the asymptotic distribution $\nu_i$ of $Y_i$ 
under $\tau_i$, respectively.
\end{Lemma}
\noindent
{\it Proof.} 
Essentially, the proof reduces to the computation for $\mu=\xi_1$.
Denote $W=11^*$ and $W_{*}=1^*1$. By Theorem 3.1, 
$$
m_k= \sum_{\pi\in\mathcal{NC}(W^k)} \prod_{V\in\pi} w(V),
$$
where $w(V) = d(V) r_{|V|}(u_V)$. Hence, $\psi_{\mu}$ takes the form
\begin{equation*}
\psi_\mu(z) = \sum_{k=1}^\infty \sum_{j_1,j_2} N_k(j_1,j_2) d_1^{j_1} d_2^{j_2} z^k,
\end{equation*}
where $j_1,j_2$ run over the set of nonnegative integers and 
$$
{\displaystyle N_k(j_1,j_2) = \sum_{\pi\in\mathcal{NC}_{j_1,j_2}(W^{k})}  
\prod_{V\in\pi} r_{|V|}(u_1)}
$$ 
where the set $\mathcal{NC}_{j_1,j_2}(W^{k})$ consists of all partitions 
$\pi\in\mathcal{NC}(W^{k})$ with the assigned dimension given by the monomial
$$
d(\pi) = \prod_{V\in\pi} d(V) = d_1^{j_1} d_2^{j_2}.
$$ 
Let us find a recurrence for the numbers $N_k(j_1,j_2)$. 
Let $\pi$ be any partition from $\mathcal{NC}_{j_1,j_2}(W^{k})$ and let 
$V_0 = \{s_1<s_2<\ldots<s_{2m}\}$ be the unique block of $\pi$ containing $s_1=1$. 
It is obvious that $s_1$ must be colored by $1$. Moreover, since $s_{i+1} - s_i$ is an odd number, 
the legs $s_1,s_3,\ldots,s_{2m-1}$ must be colored by $1$ 
and the legs $s_2,s_4,\ldots,s_{2m}$ must be colored by $1^*$. Thus, $d(V_0) = d_1^{m-1} d_2^{m}$. 
As shown in Fig. 7, there are subpartitions $\pi_1,\ldots, \pi_{2m}$ between the 
consecutive legs of $V_0$ (and after its last leg) 
such that $\pi_{2i-1} \in \mathcal{NC}(W_*^{n_{2i-1}})$ and 
$\pi_{2i} \in \mathcal{NC}(W^{n_{2i}})$ for $i\in [m]$ and some $n_1,\ldots,n_{2m}\geq 0$ 
such that $n_1 + \ldots + n_{2m} = k-m$.
From this we obtain the recurrence
$$
N_k(j_1,j_2) = \sum_{m=1}^k r_{2m} \sum_{\stackrel{j_1^{(1)}+\ldots+j_1^{(2m)} = j_1 - m +1}
{\stackrel{j_2^{(1)}+\ldots+j_2^{(2m)} = j_2 - m}{n_1 + \ldots + n_{2m} = k-m}}}  
\prod_{i=1}^m M_{n_{2i-1}}(j_1^{(2i-1)},j_2^{(2i-1)}) N_{n_{2i}}(j_1^{(2i)},j_2^{(2i)}),
$$
where 
$$
{\displaystyle M}_n(j_1,j_2) = \sum_{\pi\in\mathcal{NC}_{j_1,j_2}(W_*^{n})}\prod_{V\in\pi} r_{|V|},
$$ 
is the counterpart of ${\displaystyle N}_{n}(j_1,j_2)$ associated with the word $W_{*}^n$ and
$(r_n)$ is the sequence of free cumulants of $\nu_1$. 
We set $M_0(0,0) = N_0(0,0) = 1$ and 
$M_0(j_1,j_2) = N_0(j_1,j_2) = 0$ otherwise.
One can observe that 
$$
N_n(j_1,j_2+1)=M_n(j_1+1,j_2)
$$ 
for any $n\ge1$ and any $j_1,j_2$ since there is a natural bijection 
$$
\mathcal{NC}_{j_1,j_2+1}(W^{n})\cong \mathcal{NC}_{j_1+1,j_2}(W_*^{n})
$$ 
for any $n$ obtained by flipping the first leg (colored by $1$) 
to the last position. On the circular diagram, this bijection is natural since it 
corresponds to rotating the diagram clockwise by angle $\pi/n$. 
Hence, if we denote by $\varphi(z)$ the generating function 
without constant term of the sequence 
$(M_m(j_1,j_2))$, then 
$$
d_1\psi_\mu(z) = d_2 \varphi(z).
$$
Using this fact and the above recurrence relation, we obtain
\begin{eqnarray*}
d_1\psi_\mu 
&=&  
d_1 \sum_{k=1}^\infty  \sum_{j_1,j_2} N_k(j_1,j_2) d_1^{j_1} d_2^{j_2} z^k \\
&=&  
\sum_{k=1}^\infty \sum_{m=1}^k \sum_{j_1,j_2} r_{2m} 
\sum_{\stackrel{j_1^{(1)}+\ldots+j_1^{(2m)} = j_1}
{\stackrel{j_2^{(1)}+\ldots+j_2^{(2m)} = j_2}{n_1 + \ldots + n_{2m} = k-m}}}
\prod_{i=1}^m M_{n_{2i-1}}(j_1^{(2i-1)},j_2^{(2i-1)}) \\
& & 
\times N_{n_{2i}}(j_1^{(2i)},j_2^{(2i)}) d_1^{j_1+m} d_2^{j_2+m} z^k \\
&=&   
\sum_{m=1}^\infty  r_{2m} (zd_1d_2)^{m} \sum_{k=0}^\infty \sum_{j_1,j_2} 
\sum_{\stackrel{j_1^{(1)}+\ldots+j_1^{(2m)} = j_1}
{\stackrel{j_2^{(1)}+\ldots+j_2^{(2m)} = j_2}
{n_1 + \ldots + n_{2m} = k}}}  
\prod_{i=1}^m M_{n_{2i-1}}(j_1^{(2i-1)},j_2^{(2i-1)}) \\
& & 
\times N_{n_{2i}}(j_1^{(2i)},j_2^{(2i)}) d_1^{j_1} d_2^{j_2} z^{k} \\
&=&   
\sum_{m=1}^\infty r_{2m} \big(z d_1 d_2 (\varphi+1) (\psi_{\mu}+1)\big)^{m}  \\
&=& 
R_{{\widetilde{\nu}}}\big(z(d_1\psi_\mu+d_1)(d_1\psi_\mu+d_2)\big),
\end{eqnarray*}
where $r_n$ stands for $r_n(u_1)$. In the first equation, 
we used the fact that each segment $(1,1^*)$ of $V_0$ contributes $d_2$ whereas
each segment $(1^*,1)$ of $V_0$ contributes $d_1$. 
This proves the assertion for $\mu=\xi_1$
and $\widetilde{\nu} = \widetilde{\nu}_1$. The above formula yields
$$
R^{-1}_{\widetilde{\nu}_1}(z) = \psi^{-1}_{\xi_1}(d_1^{-1}z) (z+d_1)(z+d_2).
$$
Moreover, 
$$
\psi^{-1}_{\xi_1}(d_1^{-1} z) = \frac{z}{z+d_1} S_{\xi_1}(d_1^{-1} z) = \frac{z}{z+d_1} S_{V_{d_1}\xi_1}(z)
$$
Comparing the above equations and using 
the equation $S_{\xi_1}(z) = R^{-1}_{\xi_1}(z)/z$, we obtain
$$
S_{\xi_1}(d_{1}^{-1}z) = \frac{S_{\widetilde{\nu}_1}(z)}{z+d_2},
$$
which gives $V_{d_{1}}(\xi_1)=\widetilde{\nu_1}\boxtimes \rho_{d_2}$.
The proof for any $i\in [p]$ is similar. \hfill$\blacksquare$\\

Let us observe that the main equation for $\psi_{\mu}$ 
in the proof of Lemma 5.1 resembles that obtained by Benaych-Georges 
in [6, Lemma 3.4]. We cannot give a direct correspondence 
between these formulas since our proof is based on a different formalism, 
in which we use only free cumulants instead of their two rectangular deformations.

\begin{Lemma}
Under the assumptions of Proposition 5.1, the moment generating function $\psi_{\mu}$ of the 
asymptotic distribution $\mu$ of $BB^*$ under $\tau_1$ is the unique solution of the equation
$$
d_1\psi_\mu = R_{{\widetilde{\nu}}} \big(z (d_1\psi_\mu + d_1)(d_1\psi_\mu + d_2)\ldots(d_{1}\psi_\mu + d_{p+1}) \big)
$$
where ${\widetilde{\nu}} = {\widetilde{\nu}}_1 \boxtimes \ldots \boxtimes {\widetilde{\nu}}_p$
and $\widetilde{\nu}_j$ is defined by the even free cumulants of $\nu_j$ for $j\in [p]$.
\end{Lemma}
{\it Proof.}
Using Proposition 5.1 and Lemma 5.1, we obtain
$$
V_{d_1}(\mu)=\widetilde{\nu}
\boxtimes \rho_{d_2}\boxtimes \ldots \boxtimes \rho_{d_{p+1}},
$$
which, in terms of $S$-transforms, takes the form 
$$
S_{\mu}(d_1^{-1}z) 
= 
\frac{S_{\widetilde{\nu}}(z)}{(z+d_{2})(z+d_{3}) \ldots (z+d_{p+1})},
$$
or, equivalently,
$$
\psi^{-1}_{\mu}(d_1^{-1}z) 
=\frac{R^{-1}_{\widetilde{\nu}}(z)}{(z+d_1)(z+d_2) \ldots (z+d_{p+1})},
$$
which gives the assertion. The uniqueness of its solution follows 
from the Lagrange Inversion Theorem (see, for instance, [29]), 
which completes the proof.\hfill$\blacksquare$\vspace{12pt}

Using Lemma 5.2, we can express the moments of the distribution $\mu$ as linear combinations of 
a family of multivariate polynomials in asymptotic dimensions. For each natural $p$ (suppressed in our notation), 
we obtain the family of homogenous polynomials $\{P_{k,r}: k,r\in \mathbb{N}\}$ defined below. 
In particular, the sequence $(P_{k,1})_{k\in \mathbb{N}}$ is the sequence of 
{\it multivariate Fuss-Narayana polynomials}, which are the limit moments of products of 
independent rectangular GRM [16, 19]. Presence of polynomials $P_{k,r}$ with $r \neq 1$ 
in the formula for limit moments indicates that the considered product of Wishart type 
belongs to a different domain of attraction than the product of independent GRM. 

\begin{Definition}
{\rm By {\it generalized multivariate Fuss-Narayana polynomials} we shall 
understand polynomials of the form
$$
P_{k,r}(d_1,d_2,\ldots,d_{p+1}) = \sum_{j_1+\ldots+j_{p+1}=kp+r} 
\frac1k 
{k \choose j_1}\ldots {k \choose j_{p+1}}
d_1^{j_1} \ldots d_{p+1}^{j_{p+1}}.
$$
where $n,r\in {\mathbb N}$ and summation runs over nonnegative integers.
We use similar notation to that in [16], with variables $d_1, \ldots , d_{p+1}$, 
rather than that in [19], with $d_{0}, \ldots , d_p$.}

\end{Definition}

Let us recall the Lagrange Inversion Theorem, which says that if $f$ and $g$ are formal power series in $z$, 
for which
$$
f(z) = z g(f(z)),
$$
and $g(0)\neq0$, then
$$
[z^k] f(z) = \frac1k [\lambda^{k-1}] g^k(\lambda),
$$ 
for any $k\in \mathbb{N}$, where $[z^k] f(z)$ is the coefficient of 
$f$ standing by $z^k$ (see, for instance, [29]).

\begin{Theorem}
Moments of the asymptotic distribution $\mu$ of $BB^*$ under $\tau_1$ are given by the formula
$$
m_k = d_1^{-1}\sum_{r=1}^{k} P_{k,r}(d_1,d_2,\ldots,d_{p+1}) T_{k,r} (t_1,t_2,\ldots,t_{r})
$$
where 
$$
T_{k,r} (t_1,t_2,\ldots,t_{r}) = \sum_{i_1+\ldots+i_k=r-1} t_{i_1} t_{i_2} \ldots t_{i_k}.
$$
for any $k,r\in {\mathbb N}$, where numbers $t_k$ are coefficients of 
the $T$-transform $T_{\widetilde{\nu}}(z)=\sum_{k=0}^{\infty}t_{k}z^k$.
\end{Theorem}

\noindent
{\it Proof.} 
Using the relation $S_{\widetilde{\nu}}(z) = z^{-1}R^{-1}_{\widetilde{\nu}}(z)$,
where $R^{-1}_{\widetilde{\nu}}(z)$ is the composition inverse of $R_{\widetilde{\nu}}(z)$,
we can write the equation of Lemma 5.2 in the form
$$
d_1\psi_\mu(z) =  z (d_1\psi_\mu(z) + d_1)\ldots(d_1\psi_\mu (z)+ d_{p+1}) 
T_{{\widetilde{\nu}}}(d_1\psi_\mu(z)) .
$$
We apply the Lagrange Inversion Theorem to this equation to calculate the coefficients 
of $\psi_\mu$, which gives
\begin{eqnarray*}
d_1 m_k(\mu) \!\!&=&\!\! \frac1k [z^{k-1}] (z+d_1)^k \ldots (z+d_{p+1})^k T_{\widetilde{\nu}}^k(z) \\
\!\!
&=&
\!\! \frac1k [z^{k-1}] \sum_{j_1, \ldots , j_{p+1}=0}^k
{k\choose j_1}\ldots {k\choose j_{p+1}}  
d_1^{j_1}\ldots d_{p+1}^{j_{p+1}}  \\
&&
\sum_{r=0}^\infty \sum_{i_1+\ldots+i_k=r} 
t_{i_1} \ldots t_{i_k} z^{(p+1)k+r-j_1-\ldots -j_{p+1}}\\
\!\!
&=&
\!\! \frac1k [z^{k-1}] 
\sum_{l=0}^\infty \sum_{j_1+\ldots+j_{p+1}=l} 
{k \choose j_1}\ldots {k \choose j_{p+1}}
d_1^{j_1}\ldots d_{p+1}^{j_{p+1}}\\
&& 
\sum_{r=0}^\infty T_{k,r+1}(t_0,t_1,\ldots,t_{r}) z^{(p+1)k-l+r} \\
&=&
\!\! 
\frac1k \sum_{r=0}^{k-1} \sum_{j_1+\ldots+j_{p+1}=kp+1+r}
{k \choose j_1}\ldots {k \choose j_{p+1}} 
d_1^{j_1} \ldots d_{p+1}^{j_{p+1}}\\
&& 
T_{k,r+1}(t_0,t_1,\ldots,t_{r}) \\
&=&
\!\! \sum_{r=1}^{k} P_{k,r}(d_1,d_2,\ldots,d_{p+1}) \ T_{k,r}(t_0,t_1,\ldots,t_{r-1}),
\end{eqnarray*}
where indices $i_1, \ldots, i_n, j_1, \ldots, j_{p+1}$ are assumed to be nonnegative integers in all formulas.
This completes the proof. \hfill $\blacksquare$

\section{Examples}

In this section we present the special cases of Theorem 5.1 in which 
each measure $\nu_j$ is either the Wigner semicircle measure 
or the Marchenko-Pastur measure.

\begin{Definition} 
Let $\bm{j}=(j_1, j_2, \ldots , j_{p+1})$, where $j_1, j_2, \ldots , j_{p+1}$ are nonnegative integers. 
By $\mathcal{NC}_{\bm{j}}(W^{k})$ we denote the set of all partitions $\pi$ from $\mathcal{NC}(W^{k})$, 
for which $d(\pi) = d_1^{j_1} d_2^{j_2} \ldots d_{p+1}^{j_{p+1}}$. We define the numbers
$$
N_k(\bm{j}) = \sum_{\pi \in \mathcal{NC}_{\bm{j}}(W^{k})} \prod_{V\in\pi} r_{|V|}(u_V).
$$
for $k \ge 1$ and we set $N_0(\bm{j}) = 1$ if $\bm{j}=\bm{0}$ and $N_0(\bm{j}) = 0$ in other cases. 
Moreover, we put $N_k(\bm{j})=0$ if the set $\mathcal{NC}_{\bm{j}}(W^{k})$ is empty.
\end{Definition} 

By Theorem 3.1, the moments of the measure $\mu$ 
can be described in terms of the numbers $N_k(\bm{j})$ in the following way 
$$
m_k = \sum_{j_1,\ldots,j_{p+1}} 
N_k(\bm{j}) d_1^{j_1} d_2^{j_2} \ldots d_{p+1}^{j_{p+1}},
$$
where $j_1, \ldots , j_{p+1}$ are nonnegative integers.\\

\begin{center}
{\it Gaussian case}
\end{center}
In the first example we assume that each $\nu_i$ is the standard Wigner measure, 
thus $r_n(\nu_i)=\delta_{n,2}$ for any $i\in [p]$. 
In this case, the number $N_k(\bm{j})$ coincides with the number of all those pair partitions adapted to $W^k$ which have $j_k$ left legs colored by $k$ or 
$(k+1)^*$ for any $k\in [p+1]$, see [19].
The $R$-transform of $\widetilde{\nu}_i$ is 
simply $R_{\widetilde{\nu}_i} (z) = z$ for $i\in[p]$. 
Hence $T_{\widetilde{\nu}} (z) = 1$ and $R_{\widetilde{\nu}} (z) = z$. 
So, the equation from Lemma 5.2 takes the form
$$
d_1 \psi_\mu = z (d_1 \psi_\mu + d_1)\ldots( d_1 \psi_\mu + d_{p+1}).
$$
Note that the same equation was used in the proof of Theorem 11 in [19]. 
The only non-zero coefficient $t_n$ of 
$T_{\widetilde{\nu}}$ is $t_0 = 1$ and thus, by Theorem 5.1, the moments of $\mu$ are given by 
$$
m_k = d_1^{-1} P_{k,1}(d_1,\ldots,d_{p+1})
$$
and hence
$$
N_k(\bm{j}) = 
\frac{1}{k} 
{k \choose j_1+1}{k \choose j_2}\ldots {k \choose j_{p+1}}
$$
for $k\geq 1$ and $j_1 + j_2 + \ldots + j_{p+1} = pk$. On this occasion, let us remark that 
similar numbers were used by Edelman [10] to count chains of a certain type in the lattice $\mathcal{NC}(k)$.

\begin{Corollary}
With the above notations, $|\mathcal{NC}^2(W^k)|=F_k(p)$ for any $k,p\in {\mathbb N}$.
\end{Corollary}

\noindent
{\it Proof.} It is enough to sum all the numbers $N_k(\bm{j})$ over all $\bm{j}$ 
such that $j_1 + j_2 + \ldots + j_{p+1} = kp$, 
using Cauchy's identity (or, Vandermonde convolution). \hfill$\blacksquare$\\

\begin{center}
{\it Mixed Gaussian and Wishart case}
\end{center}
In this example we assume that the distribution $\nu_i$ is the Wigner semicircle distribution for some subset of $[p]$ and 
the standard Marchenko-Pastur (free Poisson) distribution otherwise.

Let $K$ be any subset of $[p]$ and $q$ the number of its elements. We choose that $\nu_i$ is the Wigner semicircle 
distribution (i.e. $r_n(u_i) = \delta_{n,2}$) for $i \notin K$ and $\nu_i$ is the free Poisson distribution 
(i.e. $r_n(u_i) = 1$ for any $n$) for $i\in K$. By the definition of $\widetilde{\nu}_i$, their $R$-transforms are given by
$$
R_{\widetilde{\nu}_i}(z) =
\left\{
\begin{array}{cl}
z(1-z)^{-1} & {\rm if} \;i\in K\\
z & {\rm if}\;i\notin K 
\end{array}
\right.
$$
Since $T_{\widetilde{\nu}_i}(z) = z/R^{-1}_{\widetilde{\nu}_i}(z)$ and $T_{\widetilde{\nu}}(z) = T_{\widetilde{\nu}_1}(z) \cdot \ldots \cdot T_{\widetilde{\nu}_p}(z)$, we obtain
$$
T_{\widetilde{\nu}}(z) = (1+z)^q.
$$

\begin{Proposition}
Under the above assumptions, 
$$
N_k(\bm{j}) = \frac{1}{k} 
{k \choose j_1+1}{k \choose j_2}\ldots {k \choose j_{p+1}}{kq \choose |\bm{j}|-kp},
$$
for $k\geq 1$, where $|\bm{j}|= j_1 + j_2 + \ldots + j_{p+1}$.
\end{Proposition}

\noindent
{\it Proof.} Since $T_{\widetilde{\nu}}(z) = (1+z)^q$, its coefficients are equal to 
$$
t_k = {q \choose k} 
$$
for $0\le k \le q$.
Hence, the numbers $T_{k,r}$ from Theorem 5.1 are given by
$$
T_{k,r} (t_0,t_1,\ldots,t_{r-1}) = 
\sum_{i_1+\ldots+i_k=r-1} t_{i_1} t_{i_2} \ldots t_{i_k} = 
{kq \choose r-1}
$$
and 
\begin{eqnarray*}
m_k &=& 
\sum_{r=1}^k \sum_{j_1+\ldots+j_{p+1}=kp+r} 
\frac1k
{k \choose j_1}\ldots {k \choose j_{p+1}}{kq \choose r-1}
d_1^{j_1-1} d_2^{j_2} \ldots d_{p+1}^{j_{p+1}} \\
&=& 
\sum_{j_1,\ldots,j_{p+1}\ge0} \frac1k 
{k \choose j_1+1}\ldots {k \choose j_{p+1}}{kq \choose |\bm{j}|-kp}
d_1^{j_1} d_2^{j_2} \ldots d_{p+1}^{j_{p+1}},
\end{eqnarray*}
where we assumed that 
${k \choose m}=0$ 
unless $0\le m \le k$. But the coefficient of the moment $m_k$ standing by 
$d_1^{j_1} d_2^{j_2} \ldots d_{p+1}^{j_{p+1}}$ is $N_k(\bm{j})$, which, compared with the 
last formula, gives the assertion. \hfill$\blacksquare$\\

Clearly, if $q=0$, then we get the same expression for $N_k(\bm{j})$ as in the Gaussian case. 
On the other hand, if $q=p$, then $N_k(\bm{j})$ is the cardinality of 
$\mathcal{NC}_{ \bm{j}}(W^k)$, from which the enumeration result given below follows.

\begin{Corollary}
With the above notations, $|\mathcal{NC}(W^k)| = F_k(2p)$ for any $k,p\in \mathbb{N}$.
\end{Corollary}

\begin{Example}
{\rm Consider the case where $d_i = 1$ for $i\in [p]$ and $d_{p+1} = t$ for some $t>0$. 
By Proposition 5.1 and Lemma 5.1, the measure $\mu$ is then given by
$$
\mu = \rho_1^{\;\boxtimes \,p+q-1}  \boxtimes \rho_t,
$$
hence it is the free Bessel law $\pi_{st}$ for $s = p+q$, as defined in [4]. 
By [4, Theorem 5.2], its moments are Fuss-Narayana polynomials of variable $t$. 
Indeed, by Proposition 6.1, we have
\begin{eqnarray*}
m_k &=& \sum_{j_1,\ldots,j_{p+1}\ge0} 
\frac{1}{k} 
{k\choose j_1+1}{k\choose j_2}\ldots {k\choose j_{p+1}} 
{qk \choose |\bm{j}|-pk}t^{j_{p+1}}  \\
&=& 
\sum_{j_{p+1}=1}^{k}  
{k\choose j_{p+1}}t^{j_{p+1}} 
\sum_{j_1,\ldots,j_{p}\ge0} \frac{1}{k} 
{k\choose j_1+1}{k\choose j_2} \ldots {k\choose j_{p}} {qk\choose (p+q)k - |\bm{j}|}\\
&=& 
\sum_{j_{p+1}=1}^{k} \frac{1}{k}  
{(p+q)k\choose (p+q)k+1 - j_{p+1}} 
{k\choose j_{p+1}}t^{j_{p+1}}  \\
&=& 
\sum_{j=1}^{k} \frac{1}{k} 
{(p+q)k\choose j-1}{k\choose j}t^{j}.
\end{eqnarray*}
}
\end{Example}

\begin{Example}
{\rm Assume that $d_i=t$ for all $i\in [p+1]$. Then, the measure $\mu$ is given by
$$
\mu = V_{1/t}(\rho_1^{\;\boxtimes q}  \boxtimes \rho_t^{\;\boxtimes p}).
$$ 
By Proposition 6.1, the corresponding moments $m_k$ take the form
\begin{eqnarray*}
m_k &=& \sum_{j_1,\ldots,j_{p+1}\ge0} \frac{1}{k} 
{k\choose j_1+1} {k\choose j_2}\ldots {k\choose j_{p+1}} {qk\choose |\bm{j}| - pk}t^{|\bm{j}|}  \\
&=& 
\sum_{j=pk}^{(p+q)k} \sum_{j_1+\ldots+j_{p+1}=j} \frac{1}{k} 
{k\choose j_1+1}{k\choose j_2} \ldots {k\choose j_{p+1}}{qk\choose (p+q)k - j} t^{j} \\
&=& 
\sum_{j=pk}^{(p+1)k-1} \frac{1}{k} 
{(p+1)k\choose j+1}{qk\choose  (p+q)k - j}t^{j} \\
&=& 
\sum_{j=1}^{k} \frac{1}{k} 
{(p+1)k\choose j-1}{qk\choose k-j}t^{(p+1)k-j}.
\end{eqnarray*}
In particular, if $q=1$, we obtain rescaled Fuss-Narayana polynomials of variable $t^{-1}$.
}
\end{Example}

\section{Products of dependent blocks}

Let us consider now the matrix product built from (in general, dependent) two blocks 
of the same matrix. In general, we can only find a formula for the moment generating function. 
In particular, we consider an example, in which the limit moments are rescaled Raney numbers. 

We will consider the matrix $B = X_1 X_2$, where $X_1$ and $X_2$ are square blocks of $Y\in \mathpzc{Y}$.
In other words, we assume that $p=2$ and that all asymptotic dimensions are equal to one, 
$d_1 = d_2 = d_3 = 1$.
We will denote the limit measure of $Y$ by $\nu$ and its free cumulants by
$r_k$. As before, the asymptotic distribution of $BB^*$ under $\tau_1$ will be denoted by $\mu$ 
and its moments by $m_k$, as before. We start by showing a result analogous to 
Lemma 5.1 for independent blocks. By Theorem 3.1,
$$
m_k = \sum_{\pi \in \mathcal{NC}(W^k)} \prod_{V\in\pi} r_{|V|},
$$
where $W = 122^*1^*$, but all letters here have the same label, omitted in all notations.
We are also going to use the word $W_2 = 2^*1^*12$, which is obtained from $W$ by a cyclic shift. 
Obviously, cyclic shifts of $W^k$ do not change sizes of blocks of partitions adapted to it, hence
$$
m_k = \sum_{\pi \in \mathcal{NC}(W_2^k)} \prod_{V\in\pi} r_{|V|}.
$$
\begin{Remark}
{\rm Let us notice here that the cardinality of $\mathcal{NC}^2(W^k)$ 
was considered in [15, 19] and is equal to the Fuss-Catalan number 
$F_{k}(p)$. It is natural to ask about the cardinality of $\mathcal{NC}(W^k)$. 
The study of the moments $m_k$ gives us the answer for $p = \{1,2\}$. 
For $p=1$, we obtain Edelman's result [10] that the number of 2-divisible non-crossing partitions of the set $[2k]$ 
is the Fuss-Catalan number $F_k(2)$. For $p=2$, the number of the corresponding partitions can be expressed by the 
Raney numbers, as we will demonstrate. The cases $p\ge3$ seem to go beyond the family of Raney numbers and we do not know 
their explicit form. A similar situation took place in [13].}
\end{Remark}

Let us define the following sequences of words
\begin{eqnarray*}
U_n &=& \begin{cases}
(w_1^*w_2)^m & \!\!\text{for } n=2m \\
(w_1^*w_2)^mw_1^* & \!\!\text{for } n=2m+1 \\
\end{cases} \\
U'_n &=& \begin{cases}
(w_2w_1^*)^m & \!\!\text{for } n=2m \\
(w_2w_1^*)^m w_2 & \!\!\text{for } n=2m+1 \\
\end{cases} 
\end{eqnarray*}
for $n\ge0$, where $w_1^*=1^*1$ and $w_2=22^*$. Observe that $U_n'$ is obtained from $U_n$ by 
interchanging $1$ with $2^*$ and $1^*$ with $2$. 
Therefore, the following sequence $Q_n$ can be defined in two equivalent ways:
$$
Q_n := \sum_{\pi \in \mathcal{NC}(U_n)} 
\prod_{V\in\pi} r_{|V|} = 
\sum_{\pi \in \mathcal{NC}(U'_n)} \prod_{V\in\pi} r_{|V|},
$$
for $n\ge1$, and we set $Q_0 = 1$. Obviously, $m_k = Q_{2k}$. 

We shall use the decomposition of the formal power series 
$$
\psi(z)=\sum_{n=1}^{\infty}Q_{n}z^{n}
$$ 
into the sum of
$$
\psi^{(s)}(z) = \sum_{k=1}^\infty Q_{2k} z^{2k}\;\;\;{\rm and}\;\;\;
\psi^{(a)}(z) = \sum_{k=1}^\infty Q_{2k-1} z^{2k-1},
$$
the symmetric and antisymmetric parts of $\psi$, respectively. 

\begin{Lemma}
If $\nu$ is the asymptotic distribution of $Y$, then 
$$
\psi(z) = R_{\widetilde{\nu}} \big( z (\psi(z) + 1)(\psi^{(s)}(z)+1) \big).
$$
where $\widetilde{\nu}$ is the distribution defined by the even free cumulants of $\nu$. 
\end{Lemma}

\noindent
{\it Proof.} We start by proving the following recurrence formula for the numbers $Q_n$:
$$
Q_n = \sum_{k=1}^{n} r_{2k} \sum_{2(n_1 + \ldots+ n_k) + n_{k+1} + \ldots + n_{2k} = n-k} Q_{2n_1} \ldots Q_{2n_k} Q_{n_{k+1}} \ldots Q_{n_{2k}}.
$$
for $n\ge1$. Consider $\pi\in \mathcal{NC}(U_{n})$ and its block $V_0=\{s_1 < s_2 < \ldots < s_{2k}\}$ for some $1\le k \le n$ and such that $s_1 = 1$ (the leftmost block). Observe that $s_1$ and $s_2$ must be colored by $1^*$ and $1$, respectively, 
and that the word between them is $W^{m_1}$ for some $m_1\ge0$. Hence, the contribution from 
all possible `subpartitions' of $\pi$ on the set $\{2,3,\ldots,s_2-1\}$ to the moment $Q_n$ is equal to $Q_{2m_1}$. 

Now, $s_3$ may be colored by $1^*$ or $2$. In both cases the word between $s_2$ and $s_3$ is $U'_{m_{2}}$ for some $m_2\ge0$ and the contribution from all possible `subpartitions' of $\pi$ on the set $\{s_2+1,\ldots,s_3-1\}$ is $Q_{m_2}$. Similarly, we can see that the word between $s_3$ and $s_4$ must be $W^{m_3}$ or $W_2^{m_3}$ and their contribution is $Q_{2m_3}$, and so on. Changing the sequence of indices by: $m_{2i-1}$ to $n_i$ and $m_{2i}$ to $n_{i+k}$ for $i\in[k]$ and counting the letters of the word $U_n$ with respect to the above decomposition, we arrive at
$$
2(n_1 + \ldots+ n_k) + n_{k+1} + \ldots + n_{2k} = n-k,
$$
which proves the desired recurrence formula for $Q_n$. Using it, we obtain
\begin{eqnarray*}
\psi(z) &=& \sum_{n=1}^\infty Q_n z^n \\
%&=& \sum_{n=1}^\infty \sum_{k=1}^{n} r_{2k} \sum_{2n_1 + \ldots+ 2n_k + n_{k+1} + \ldots + n_{2k} = n-k} Q_{2n_1} \ldots Q_{2n_k} Q_{n_{k+1}} \ldots Q_{n_{2k}} z^n \\
&=& \sum_{k=1}^\infty r_{2k} z^k \sum_{n=0}^{\infty} \sum_{2n_1 + \ldots+ 2n_k + n_{k+1} + \ldots + n_{2k} = n} Q_{2n_1} \ldots Q_{2n_k} Q_{n_{k+1}} \ldots Q_{n_{2k}} z^n \\
&=& \sum_{k=1}^\infty r_{2k} z^k (\psi(z) + 1)^k 
(\psi^{(s)}(z) + 1)^k = R_{\widetilde{\nu}} \big( z (\psi(z) + 1)(\psi^{(s)}(z)+1) \big),
\end{eqnarray*}
which ends the proof. \hfill$\blacksquare$

\begin{Example}
{\rm Assume that $X_1 = X_2$ and that $\nu$ is the standard Wigner semicircle law. 
Then $R_{\widetilde{\nu}}(z) = z$ and Lemma 7.1 takes the form
$$
\psi^{(a)}(z) + \psi^{(s)}(z) = z (\psi^{(a)}(z) + \psi^{(s)}(z) + 1)(\psi^{(s)}(z)+1),
$$
which gives 
$$
\begin{cases}
\psi^{(s)}(z) = z \psi^{(a)}(z) (\psi^{(s)}(z)+1)\\ 
\psi^{(a)}(z) = z (\psi^{(s)}(z)+1)^2
\end{cases}
$$
and thus
$$
\psi^{(s)}(z) = z^2 (\psi^{(s)}(z)+1)^3.
$$
Moreover, $\psi^{(s)}(z) = \psi_\mu(z^2)$ since $m_k = Q_{2k}$ for any $k\in \mathbb{N}$. 
Therefore, 
$$
\psi_\mu(z) = z (\psi_\mu(z) + 1)^3.
$$
This means that $m_n=F_n(3)$ for any $n\in \mathbb{N}$, 
which agrees with the more general result in [3], where it was shown that the $p$th power of a GRM $X$ behaves asymptotically as the product of 
$p$ independent copies of $X$.}
\end{Example}

Let us recall that the {\it Raney numbers} are given by the formula
$$
R_n(p,r) = \frac{r}{(p+1)n+r} {(p+1)n+r\choose n},
$$
for $p,r\in\mathbb{R}_+$, where
$$
{x \choose n}=x(x-1)\ldots(x-n+1)/n!
$$
for any $x\in\mathbb{R}$. In particular, if $p=r=1$, we obtain Catalan numbers.

\begin{Example}
{\rm 
Assume that $\nu$ is the free Poisson distribution. Then 
$r_n = 1$ for $n\ge1$ and hence $R_{\widetilde{\nu}}(z) = z(1-z)^{-1}$.
The equation of Lemma 7.1 takes the form
$$
\psi(z) = \frac{z(\psi(z)+1)(\psi^{(s)}(z)+1)}{1-z(\psi(z)+1)(\psi^{(s)}(z)+1)},
$$
which, when multiplied by the denominator of the RHS, gives
$$
\psi(z) = z (\psi(z) + 1)^2 (\psi^{(s)}(z) + 1).
$$
Decomposition into the symmetric and antisymmetric parts gives
$$
\begin{cases}
\psi^{(s)}(z) = 2 z \psi^{(a)}(z) (\psi^{(s)}(z)+1)^2 \\ 
\psi^{(a)}(z) = z \big( (\psi^{(s)}(z)+1)^3 + (\psi^{(a)}(z))^2 (\psi^{(s)}(z)+1) \big),
\end{cases}
$$
which leads to the equation
$$
(\psi^{(s)}(z)+1)^2 = 1 + 4z^2 (\psi^{(s)}(z)+1)^6, 
$$
or
$$
M_\mu(z)^2 = 1 + 4z M_\mu(z)^6,
$$
where $M_\mu(z) = \sum_{n=0} m_n z^n$.
Applying the Lagrange Inversion Theorem to the last formula, we get 
\begin{eqnarray*}
[z^{n}] M_\mu(z)^2 &=&
\frac1n [\lambda^{n-1}] 4^n (\lambda+1)^{3n}\\
&=&
 \frac{4^n}{n}  
{3n \choose n-1}= 
\frac{4^n}{3n+1}  
{3n+1 \choose n}.
\end{eqnarray*}
On the other hand, the coefficient $[z^{n}] M_\mu(z)^2$ is the convolution of the moment
sequence $(m_k)$ with itself, thus
$$
\sum_{k=0}^n m_k m_{n-k}  = \frac{4^n}{3n+1}  
{3n+1 \choose n}.
$$
The last formula can be viewed as a recurrence relation 
with indeterminate $m_n$. Obviously, this equation has a unique solution and by the well known property of Raney numbers
$$
\sum_{k=0}^n R_k(p,r) R_{n-k}(p,s) = R_n(p,r+s)
$$
we get
$$
m_n = 4^n R_n(2,1/2) = \frac{4^n}{6n+1}  
{3n+\frac12 \choose n}.
$$
}
\end{Example}

From this example one can immediately obtain a new enumeration result, similar to that
in Corollary 6.2 for $p=2$. When comparing both results, one has to remember 
that the class $\mathcal{NC}(W^{n})$ is now different than that used in Corollary 6.2 
since all letters in $W=122^*1^*$ have the same label, which was not the case before.
The corresponding twice longer word $\widetilde{W}$ that is suitable for pair partitions 
can be identified with $W_{0}=122^*44^*22^*1^*$, as explained in Section 4.

\begin{Corollary}
With the above notations, 
$$
|\mathcal{NC}(W^n)| = |\mathcal{NC}^2(\widetilde{W}^n)|=4^n R_{n}(2,1/2)$$ 
for any $n\in \mathbb{N}$.
\end{Corollary}

\noindent\\[10pt]
{\bf Acknowledgement}\\
We would like to thank Wojtek M{\l}otkowski for an observation that led to Remark 5.2.

\end{document}